# MATRIX NORMS AND RAPID MIXING FOR SPIN SYSTEMS

By Martin Dyer, Leslie Ann Goldberg and Mark Jerrum

*University of Leeds, University of Liverpool and University of London*

We give a systematic development of the application of matrix norms to rapid mixing in spin systems. We show that rapid mixing of both random update Glauber dynamics and systematic scan Glauber dynamics occurs if *any* matrix norm of the associated dependency matrix is less than 1. We give improved analysis for the case in which the diagonal of the dependency matrix is $\mathbf{0}$ (as in heat bath dynamics). We apply the matrix norm methods to random update and systematic scan Glauber dynamics for coloring various classes of graphs. We give a general method for estimating a norm of a symmetric nonregular matrix. This leads to improved mixing times for any class of graphs which is hereditary and sufficiently sparse including several classes of degree-bounded graphs such as nonregular graphs, trees, planar graphs and graphs with given tree-width and genus.

**1. Introduction.** A *spin system* consists of a finite set of *sites* and a finite set of *spins*. A *configuration* is an assignment of a spin to each site. Sites interact locally, and these interactions specify the relative likelihood of possible (local) subconfigurations. Taken together, these give a well-defined probability distribution $\pi$ on the set of configurations.

*Glauber dynamics* is a Markov chain whose states are configurations. In the transitions of the Markov chain, the spins are updated one at a time. The Markov chain converges to the *stationary distribution* $\pi$. During each transition of *random update Glauber dynamics*, a site is chosen *uniformly at random* and a new spin is chosen from an appropriate probability distribution (based on the local subconfiguration around the chosen site). During a transition of *systematic scan Glauber dynamics*, the sites are updated in a (deterministic) systematic order, one after another. Again, the updates are from an appropriate probability distribution based on the local subconfiguration.

It is well known that the mixing times of random update Glauber dynamics and systematic scan Glauber dynamics can be bounded in terms of the influences of sites on each other. A *dependency matrix* for a spin system with $n$ sites is an $n \times n$ matrix $R$ in which $R_{i,j}$ is an upper bound on the *influence* (defined below) of site $i$ on site $j$.

An easy application of the *path coupling* method of Bubley and Dyer shows that if the $L_\infty$ norm of $R$ (which is its maximum row sum and is written $\|R\|_\infty$) is less than 1 then random update Glauber dynamics is rapidly mixing. The same is true if the $L_1$ norm (the maximum column sum of $R$, written $\|R\|_1$) is less than 1. The latter condition is known as the *Dobrushin condition*. Dobrushin [11] showed that if $\|R\|_1 < 1$, then the corresponding countable spin system has a *unique Gibbs measure*. As we now know (see Weitz [39]), there is a very close connection between rapid mixing of Glauber dynamics for finite spin systems and uniqueness of Gibbs measure for the corresponding countable systems. Dobrushin and Shlosman [12] were the first to establish uniqueness when $\|R\|_\infty < 1$. Their analysis extends to *block dynamics* but we will stick to Glauber dynamics in this paper. For an extension of some of our ideas to block dynamics, see [30].

The Dobrushin condition $\|R\|_1 < 1$ implies that systematic scan is rapidly mixing. A proof follows easily from the account of Dobrushin uniqueness in Simon's book [35], some of which is derived from the account of Föllmer [19]. In [14], we showed that $\|R\|_\infty < 1$ also implies rapid mixing of systematic scan Glauber dynamics. [14], Section 3.5, notes that it is possible to prove rapid mixing by observing a contraction in other norms besides the $L_1$ norm and the $L_\infty$ norm. This idea was developed by Hayes [22], who showed that rapid mixing occurs when the spectral norm $\|R\|_2$ is less than one. For symmetric matrices, the spectral norm is equal to the largest eigenvalue of $R$, $\lambda(R)$. So, for symmetric matrices, [22] gives rapid mixing when $\lambda(R) < 1$. In general, $\|R_2\|/\lambda(R)$ can be arbitrarily large, see Section 2.1.

In this paper, we give a systematic development of the application of matrix norms to rapid mixing. We first show that rapid mixing of random update Glauber dynamics occurs if *any* matrix norm is less than 1. Formally, we prove the following, where $J_n$ is the norm of the all 1's matrix. All definitions are given in Section 2.

LEMMA 1. *Let $R$ be a dependency matrix for a spin system, and let $\|\cdot\|$ be any matrix norm such that $\|R\| \le \mu < 1$. Then the mixing time of random update Glauber dynamics is bounded by*

$$\hat{\tau}_r(\varepsilon) \sim n(1-\mu)^{-1} \ln((1-\mu)^{-1} J_n/\varepsilon).$$

We prove a similar result for systematic scan Glauber dynamics.

Lemma 2. *Let $R$ be a dependency matrix for a spin system and $\|\cdot\|$ any matrix norm such that $\|R\| \le \mu < 1$. Then the mixing time of systematic scan Glauber dynamics is bounded by*

$$\hat{\tau}_S(\varepsilon) \sim (1-\mu)^{-1}\ln((1-\mu)^{-1}J_n/\varepsilon).$$

The chief benefit of the new lemmas is that they can be used to show rapid mixing whenever the dependency matrix has *any* norm which is less than 1, even if the norms which are mentioned in previous theorems are not less than 1. Section 2.3 gives an example of a spin system for which Lemmas 1 and 2 can be used to prove rapid mixing, while previous theorems are inapplicable. The point of the lemmas is that rapid mixing occurs whenever *any* matrix norm is bounded—specific properties of the norm are not relevant.

Section 3.1 uses path coupling to prove Lemmas 1 and 2. Despite historical differences, the path-coupling approach is essentially equivalent to Dobrushin uniqueness. To demonstrate the relationship between the approaches, we again prove the same lemmas using Dobrushin uniqueness in Section 3.2. We also give an improved analysis for the case in which the diagonal of $R$ is $\mathbf{0}$, which is the case for *heat bath* dynamics. We prove the following.

Lemma 3. *Let $R$ be symmetric with zero diagonal and $\|R\|_2 = \lambda(R) = \lambda < 1$. Then the mixing time of systematic scan is at most*

$$\hat{\tau}_S(\varepsilon) \sim (1-\tfrac{1}{2}\lambda)(1-\lambda)^{-1}\ln((1-\lambda)^{-1}n/\varepsilon).$$

An interesting observation is that when $\lambda(R)$ is close to 1, the number of Glauber steps given in the upper bound from Lemma 3 is close to half the number that we get in our best estimate for random update Glauber dynamics (see Remark 6)—perhaps this can be interpreted as weak evidence in support of the conjecture that systematic scan mixes faster than random update for Glauber dynamics.

1.1. *Applications.* The study of spin systems originates in statistical physics. Configurations in spin systems are used to model configurations in physical systems involving interacting particles. Rapid mixing is important for two reasons.

(i) When Glauber dynamics is rapidly mixing, it can be used for sampling. Typically, we are interested in sampling configurations to learn about the equilibrium distribution. In particular, we are often interested in estimating the so-called *partition function* of the system. If Glauber dynamics is rapidly mixing, then a short simulation (of feasible length) yields a sample distribution which is close to the equilibrium distribution. Otherwise, Glauber dynamics is an inappropriate means of producing samples.

(ii) Rapid mixing of Glauber dynamics is strongly associated with qualitative properties such as uniqueness of the infinite-volume Gibbs measure. Infinite systems are beyond the scope of this paper, but it is interesting to note that there are rigorous proofs that rapid mixing of Glauber dynamics on finite systems often coincides with uniqueness, which is a qualitative property on infinite systems—the property of having one, rather than many, qualitative equilibria. See [27, 39] for more details about this fascinating connection.

In computer science, rapid mixing has important applications to the computational complexity of *counting* problems and their relatives. While exact counting seems intractable in most cases, efficient *sampling* usually implies the possibility of efficient *approximate* counting [25]. In this area, considerable attention has been paid to problems which are essentially spin systems, for example, colorings and independent sets in graphs [31]. Here the specific dynamics is not important, only that it has polynomial mixing time. However, it is generally the case that if any dynamics mixes rapidly then so will the Glauber dynamics. This can usually be established using *Markov chain comparison* techniques [9, 17, 32]. Therefore, the Glauber dynamics still retains a central importance.

Traditionally, rigorous analysis has focused on the mixing properties of random update Glauber dynamics, which is easier to analyze than systematic scan Glauber dynamics. (See [1, 8, 16] for a discussion of some notable exceptions.) However, experimental work is often carried out using systematic scan strategies. Thus, it is important to understand the mixing time of systematic scan Glauber dynamics. The observation that the Dobrushin condition implies that systematic scan is rapidly mixing (which is an observation of Sokal) was an important breakthrough. This was extended in [14] which showed that the Dobrushin–Shlosman condition (bounding the $L_\infty$ norm) also implies rapid mixing of systematic scan Glauber dynamics. Dyer, Goldberg and Jerrum [14] gave an application to sampling proper colorings of an arbitrary degree bounded graph. This is an important application in computer science because colorings are used to model many combinatorial structures such as assignments and timetables.

Hayes [22] gives applications of conditions of the Dobrushin type to various related problems on graphs, using the norm $\|\cdot\|_2$. In [14], we observed that the dependency matrix for the Glauber dynamics on graph colorings can be bounded by a multiple of the adjacency matrix of the graph. This was applied to analyzing the systematic scan dynamics for coloring near-regular graphs, and hence to regular graphs. Hayes extends the observation of [14] to the Glauber dynamics for the Ising and hard core models. He applies these observations with a new estimate of the largest eigenvalue of the adjacency matrix of a planar graph, obtaining an improved estimate of the mixing time of these chains on planar graphs with given maximum degree. He also applies them to bounded-degree trees, using an eigenvalue estimate due to Stevanović [36], for which he provides a different proof. He extends these results to the systematic scan chain for each problem, using ideas taken from [14].

In Section 4, we apply the matrix norm methods developed here to the random update Glauber dynamics and systematic scan dynamics for coloring various classes of graphs. We give a general method for estimating the norm $\|\cdot\|_2$ of a symmetric nonnegative matrix $R$. Our method is again based on matrix norms. We show that there exists a "decomposition" $R = B + B^{\mathrm{T}}$, for some matrix $B$, where $\|B\|_1$, $\|B\|_\infty$ can be bounded in terms of $\|R\|_1$ and the *maximum density* of $R$. The bounds on $\|B\|_1$, $\|B\|_\infty$ can then be combined to bound $\|R\|_2$. In particular, our methods allow us to give a common generalization of results of Hayes [22], Stevanović [36] and others for the maximum eigenvalue of certain graphs. In most cases, we are also able to strengthen the previous results. In particular, Corollaries 49(i) and 49(ii) improve the results of Stevanović and Zhang and Corollary 49(iv) improves a result of Hayes. Theorem 51 gives new rapid-mixing results for sparse hereditary graph classes.

Using this, we obtain whole classes of graphs for which we did not have rapid mixing results which improved those on arbitrary degree-bounded graphs, but now we do. These results are summarized in Corollary 52. Part (i) gives mixing time bounds for all connected graphs when $q$, the number of spins, is equal to twice the degree, $\Delta$. The $q = 2\Delta$ boundary case is important and well studied. Part (i) improves the mixing time bound given in Theorem 5 of [14] by a factor of $n$. Part (ii) gives mixing time bounds for graphs with bounded tree-width. These extend results by Martinelli, Sinclair and Weitz [28] which show rapid mixing for trees, but not for graphs with higher treewidth (trees are graphs with treewidth 1). Part (iii) gives mixing-time bounds for planar graphs. These improve the results of Hayes [22] which do not apply unless $q$ is increased by a fixed multiple of $\Psi$. The goal is to get rapid-mixing results for $q$ as small as possible. For trees, it is known that $q = \Delta + 3$ suffices, and it is an open question how small $q$ can be as a function of $\Delta$ for these other graph classes. Part (iv) improves our planar graphs results by extending them to general graphs with bounded genus, rather than just to planar graphs. Prior to our work, rapid mixing was known only for $q \geq 11\Delta/6$ [38].

**2. Preliminaries.** Let $[n] = \{1, 2, \ldots, n\}$, $\mathbb{N} = \{1, 2, 3, \ldots\}$, and $\mathbb{N}_0 = \mathbb{N} \cup \{0\}$. We use $\mathbb{Z}, \mathbb{R}$ for the integers and reals, and $\mathbb{R}_+$ for the nonnegative reals. Let $|c|$ denote the absolute value of $c$.

2.1. *Matrix norms.* Let $\mathbb{M}_{mn} = \mathbb{R}^{m \times n}$ be the set of real $m \times n$ matrices. We denote $\mathbb{M}_{nn}$, the set of *square* matrices, by $\mathbb{M}_n$. The set of *nonnegative* matrices will be denoted by $\mathbb{M}^+_{mn}$, and the set of square nonnegative matrices by $\mathbb{M}^+_n$. We will write $\mathbf{0}$ for the zero $m \times n$ matrix and $I$ for the $n \times n$ identity matrix. The dimensions of these matrices can usually be inferred from the context, but where ambiguity is possible or emphasis required, we will write $\mathbf{0}_{m,n}$, $I_n$, etc. Whether vectors are row or column will be determined either by context or explicit statement. The $i$th component of a vector $v$ will be written both as $v_i$ and $v(i)$, whichever is more convenient. If $R$ is a matrix and $v$ a vector, $Rv(i)$ will mean $(Rv)_i$. We will use $\mathbf{J}$ for the $n \times n$ matrix of 1's, $\mathbf{1}$ for the column $n$-vector of 1's, and $\mathbf{1}^{\mathrm{T}}$ for the row $n$-vector of 1's. Again, the dimensions can be inferred from context.

A *matrix norm* (see [23]) is a function $\|\cdot\| : \mathbb{M}_{mn} \to \mathbb{R}_+$ for each $m, n \in \mathbb{N}$ such that:

(i) $\|R\| = 0$ and $R \in \mathbb{M}_{mn}$ if and only if $R = \mathbf{0} \in \mathbb{M}_{mn}$;
(ii) $\|\mu R\| = |\mu|\|R\|$ for all $\mu \in \mathbb{R}$ and $R \in \mathbb{M}_{mn}$;
(iii) $\|R + S\| \leq \|R\| + \|S\|$ for all $R, S \in \mathbb{M}_{mn}$;
(iv) $\|RS\| \leq \|R\|\|S\|$ for all $R \in \mathbb{M}_{mk}$, $S \in \mathbb{M}_{kn}$ $(k \in \mathbb{N})$.

Note that property (2.1) (*submultiplicativity*) is sometimes not required for a matrix norm, but we will require it here. The condition that $\|\cdot\|$ be defined for all $m, n$ is, in fact, a mild requirement. Suppose $\|\cdot\|$ is initially defined only on $\mathbb{M}_n$, for any large enough $n$, then we can define $\|R\|$ for $R \in \mathbb{M}_{mk}$ $(m, k \in [n])$ by "embedding" $R$ in $\mathbb{M}_n$, that is,

$$\|R\| \stackrel{\text{def}}{=} \left\| \left[ \begin{array}{c|c} R & \mathbf{0}_{m,n-k} \\ \hline \mathbf{0}_{n-m,k} & \mathbf{0}_{n-m,n-k} \end{array} \right] \right\|.$$

It is straightforward to check that this definition gives the required properties. For many matrix norms, this embedding norm coincides with the actual norm for all $m, k \in [n]$.

Examples of matrix norms are *operator norms*, defined by $\|R\| = \max_{x \neq 0} \|Rx\|/\|x\|$ for any *vector norm* $\|x\|$ defined on $\mathbb{R}^n$ for all $n \in \mathbb{N}$. Observe that we denote a matrix norm by $\|\cdot\|$ and a vector norm by $\|\cdot\|$. Since vector norms occur only in this section, this should not cause confusion. In fact, their meanings will also be very close, as we discuss below.

For any operator norm, we clearly have $\|I\| = 1$. The norms $\|\cdot\|_1$, $\|\cdot\|_2$ and $\|\cdot\|_\infty$, are important examples, derived from the corresponding vector norms. The norm $\|R\|_1$ is the maximum column sum of $R$, $\|R\|_\infty$ is the maximum row sum, and the *spectral* norm $\|R\|_2 = \sqrt{\lambda}$, where $\lambda$ is the largest eigenvalue of $R^{\mathrm{T}}R$. (See [23], pages 294–295, but observe that $|||\cdot|||$ is used for what we denote here by $\|\cdot\|$.) The *Frobenius norm* $\|R\|_F = \sqrt{\sum_{i,j} R_{ij}^2}$ (see [23], page 291) is an example of a matrix norm which is not an operator norm. Note that $\|I\| = \sqrt{n}$ for the Frobenius norm, so it cannot be defined as an operator norm.

New matrix norms can also be created easily from existing ones. If $W_n \in \mathbb{M}_n$ is a fixed nonsingular matrix for each $n$, then $\|\cdot\|_W = \|W_m(\cdot)W_n^{-1}\|$ is a matrix norm. (See [23] page 296.) Note that $\|\cdot\|_W$ is an operator norm whenever $\|\cdot\|$ is, since it is induced by the vector norm $\|W_m \cdot \|$.

The following relate matrix norms to absolute values and corresponding vector norms.

LEMMA 4. *Suppose $c \in \mathbb{R}$. Let $\|\cdot\|$ be a matrix norm on $1 \times 1$ matrices. Then $|c| \le \|c\|$.*

PROOF. This follows from the axioms for a matrix norm. First, $\|c\| = \|c \times 1\| = |c|\|1\|$ by (ii). Also, $\|c\| = \|c \times 1\| \le \|c\|\|1\|$ by (iv). Finally, $\|1\| \ne 0$ by (i). $\square$

LEMMA 5. *Suppose $x$ is a column vector, $\|\cdot\|$ a vector norm and $\|\cdot\|$ the corresponding operator norm. Then $\|x\| = \|1\|\|x\|$.*

PROOF. Let $x$ be a length-$\ell$ column vector. $\|x\|$ is the vector norm applied to $x$, $\|1\|$ is the same norm applied to the length-1 column vector containing a single 1. $\|x\|$ is the operator norm applied to the $\ell \times 1$ matrix containing the single column $x$. Then $\|x\| = \max_{\alpha \ne 0} \|x\alpha\|/\|\alpha\|$ where $\alpha$ is a nonnegative real number. Pulling constants out of the vector norm, $\max_{\alpha \ne 0} \|x\alpha\|/\|\alpha\| = \|x\|/\|1\|$. $\square$

The *dual* (or *adjoint* [23], page 309) norm $\|\cdot\|^*$ of a matrix norm $\|\cdot\|$ will be defined by $\|R\|^* = \|R^{\mathrm{T}}\|$. Thus, $\|\cdot\|_1$ and $\|\cdot\|_\infty$ are dual, and $\|\cdot\|_2$ is self-dual. Note that, for any column vector $x$, $\|x^{\mathrm{T}}\| = \|x\|^*$ so, for example, $\|x^{\mathrm{T}}\|_1 = \|x\|_\infty$. Clearly, any matrix norm $\|\cdot\|$ induces a vector norm $\|\cdot\|$ on column vectors. Then the dual matrix norm, as defined here, is closely related to the *dual vector norm*, which is defined by

$$\|x\|^* = \max_{y \ne 0} \frac{|x^{\mathrm{T}} y|}{\|y\|}.$$

LEMMA 6. *Suppose $x$ is a column vector, $\|\cdot\|$ a vector norm, and $\|\cdot\|$ the corresponding operator norm. Then $\|x\|^* = \|1\|^* \, \|x\|^*$.*

PROOF. By definition, $\|1\|^* = \max_{\alpha \ne 0} |\alpha|/\|\alpha\| = 1/\|1\|$, after pulling out constants, and

$$\|1\|\|x\|^* = \|1\| \max_{y \ne 0} \frac{|x^{\mathrm{T}} y|}{\|y\|} = \max_{y \ne 0} \frac{|x^{\mathrm{T}} y|\|1\|}{\|y\|} = \max_{y \ne 0} \frac{\|x^{\mathrm{T}} y\|}{\|y\|} = \|x^{\mathrm{T}}\| = \|x\|^*. \quad \square$$

With any matrix $R = (R_{ij}) \in \mathbb{M}_n$ we can associate a weighted digraph $G(R)$ with vertex set $[n]$, edge set $E = \{(i,j) \in [n]^2 : R_{ij} \neq 0\}$, and $(i,j) \in E$ has weight $R_{ij}$. The (zero-one) adjacency matrix of $G(R)$ will be denoted by $A(R)$. If $G(R)$ is labeled so that each component has consecutive numbers, then $R$ is *block diagonal* and the (principal) *blocks* correspond to the components of $G(R)$. A block is *irreducible* if the corresponding component of $G(R)$ is *strongly connected*. Note, in particular, that $R$ is irreducible if $R > \mathbf{0}$. If $R$ is *symmetric*, $G(R)$ is an *undirected* graph and $R$ is irreducible when $G(R)$ is connected. For $i, j \in V$, $\mathrm{d}(i,j)$ will denote the number of edges in a shortest directed path from $i$ to $j$. If there is no such path, $\mathrm{d}(i,j) = \infty$. The diameter of $G$, $D(G) = \max_{i,j \in V} \mathrm{d}(i,j)$. Thus, $G$ is strongly connected when $D(G) < \infty$.

For $R \in \mathbb{M}_n^+$, let $\lambda(R)$ denote the largest eigenvalue (the *spectral radius*). We know that $\lambda(R) \in \mathbb{R}_+$ from Perron–Frobenius theory [34], Chapter 1. We use the following facts about $\lambda(R)$. The first is a restatement of [34], Theorem 1.6, a version of the Perron–Frobenius theorem.

LEMMA 7. *If $R \in \mathbb{M}_n^+$ is irreducible, there exists a row vector $w > \mathbf{0}$ satisfying $wR \leq \mu w$ if and only if $\mu \geq \lambda(R)$. If $\mu = \lambda(R)$, then $w$ is the unique left eigenvector of $R$ for the eigenvalue $\lambda$.*

LEMMA 8. *If $R \in \mathbb{M}_n^+$ has blocks $R_1, R_2, \ldots, R_k$, then $\lambda(R) = \max_{1 \leq i \leq k} \lambda(R_i)$.*

LEMMA 9 (See [34], Theorem 1.1). *If $R, R' \in \mathbb{M}_n^+$ and $R \leq R'$, then $\lambda(R) \leq \lambda(R')$.*

$\lambda(\cdot)$ is not a matrix norm. For example,

$$\lambda \begin{pmatrix} 0 & 1 \\ 0 & 0 \end{pmatrix} = 0$$

so axiom (i) in the definition of a matrix norm is violated by $\lambda(\cdot)$. Nevertheless, $\lambda(R)$ is a lower bound on the value of any norm of $R$.

LEMMA 10 (See [23], Theorem 5.6.9). *If $R \in \mathbb{M}_n^+$, then $\lambda(R) \leq \|R\|$ for any matrix norm $\|\cdot\|$.*

Furthermore, for every $R \in \mathbb{M}_n^+$ there is a norm $\|\cdot\|$, depending on $R$, such that the value of this norm coincides with $\lambda(\cdot)$ when evaluated at $R$.

LEMMA 11. *For any irreducible $R \in \mathbb{M}_n^+$, there exists a matrix norm $\|\cdot\|$ such that $\lambda(R) = \|R\|$.*

PROOF. Let $w > \mathbf{0}$ be a left eigenvector for $\lambda = \lambda(R)$, and let $W = \operatorname{diag}(w) \in \mathbb{M}_n^+$. Then $\|\cdot\|_w = \|W(\,\cdot\,)W^{-1}\|_1$ is the required norm, since $\|R\|_w = \|WRW^{-1}\|_1 = \|wRW^{-1}\|_1 = \lambda\|wW^{-1}\|_1 = \lambda\|\mathbf{1}^{\mathrm{T}}\|_1 = \lambda\|\mathbf{1}\|_\infty = \lambda$. $\square$

The norm $\|\cdot\|_w$ defined in the proof of Lemma 11 is the *minimum* matrix norm for $R$, but this norm is clearly dependent upon $R$ since $w$ is.

The *numerical radius* [23] of $R \in \mathbb{M}_n^+$ is defined as $\nu(R) = \max\{x^{\mathrm{T}}Rx : x^T x = 1\}$. $\nu(\cdot)$ is not submultiplicative since

$$\nu\begin{pmatrix} 0 & 1 \\ 0 & 0 \end{pmatrix} = \nu\begin{pmatrix} 0 & 0 \\ 1 & 0 \end{pmatrix} = \tfrac{1}{2},$$

but applying $\nu$ to the product gives

$$\nu\begin{pmatrix} 1 & 0 \\ 0 & 0 \end{pmatrix} = 1.$$

Thus, $\nu(\cdot)$ is not a matrix norm in our sense. Nevertheless, $\nu(R)$ provides a lower bound on the norm $\|R\|_2$.

LEMMA 12. *$\lambda(R) \le \nu(R) \le \|R\|_2$, with equality throughout if $R$ is symmetric.*

PROOF. Let $w$, with $\|w\|_2 = 1$, be an eigenvector for $\lambda = \lambda(R)$. Then $\nu(R) \ge w^{\mathrm{T}}Rw = \lambda w^{\mathrm{T}}w = \lambda(R)$. Also, $\nu(R) = x^{\mathrm{T}}Rx \le \|R\|_2$ for some $x$ with $\|x\|_2 = 1$, and $x^{\mathrm{T}}Rx = \|x^{\mathrm{T}}Rx\|_2 \le \|R\|_2$ since $\|\cdot\|_2$ is submultiplicative. If $R$ is a symmetric matrix, then $R = Q^{\mathrm{T}}\Lambda Q$, for $Q$ orthonormal and $\Lambda$ a diagonal matrix of eigenvalues. Then $\|R\|_2^2 = \nu(R^{\mathrm{T}}R) = \nu(\Lambda^2) = \lambda(R)^2$. $\square$

Thus, when $R$ is symmetric, we have $\lambda(R) = \|R\|_2$, and hence $\|\cdot\|_2$ is the minimum matrix norm, uniformly for all symmetric $R$. However, when $R$ is not symmetric, $\|R\|_2/\lambda(R)$ can be arbitrarily large, even though $\mathbf{0} < R < \mathbf{J}$. Consider, for example,

$$R = \begin{bmatrix} \varepsilon & 1 - 2\varepsilon \\ \varepsilon & \varepsilon \end{bmatrix},$$

for any $0 < \varepsilon < \frac{1}{2}$. Then $\lambda(R) < \sqrt{\varepsilon} + \varepsilon$, and $\|R\|_2 > 1 - 2\varepsilon$, so $\lim_{\varepsilon\to 0} \|R\|_2/\lambda(R) = \infty$. Also, $\|\cdot\|_2$ is not necessarily the minimum norm for asymmetric $R$. We always have $\|R\|_2 \le \sqrt{n}\|R\|_1$ ([23], page 314), but this bound can almost be achieved for $\mathbf{0} < R < \mathbf{J}$. Consider

$$R = \begin{bmatrix} 1 - n\varepsilon & 1 - n\varepsilon & \dots & 1 - n\varepsilon & 1 - n\varepsilon \\ \varepsilon & \varepsilon & \dots & \varepsilon & \varepsilon \\ \vdots & \vdots & & \vdots & \vdots \\ \varepsilon & \varepsilon & \dots & \varepsilon & \varepsilon \end{bmatrix},$$

for any $0 < \varepsilon < \frac{1}{n}$. Then $\|R\|_1 = 1 - \varepsilon$, but $\|R\|_2 > (1 - n\varepsilon)\sqrt{n}$, so $\lim_{\varepsilon \to 0} \|R\|_2 / \|R\|_1 = \sqrt{n}$. On the other hand, $\|\cdot\|_2$ does have the following minimality property.

LEMMA 13. *For any matrix norm* $\|\cdot\|$, $\|R\|_2 \le \sqrt{\|R\|\|R\|^*}$.

PROOF. $\|R\|_2^2 = \lambda(R^{\mathrm{T}}R) \le \|R^{\mathrm{T}}R\| \le \|R^{\mathrm{T}}\|\|R\| = \|R\|\|R\|^*$, using Lemmas 10 and 12. □

For a matrix norm $\|\cdot\|$, the quantities $J_n = \|\mathbf{J}\|$, for $\mathbf{J} \in \mathbb{M}_n$ and $C_n = \|\mathbf{1}\|\|\mathbf{1}\|^*$, will be used below. We collect some of its properties here. In particular, $J_n = n$ for $\|\cdot\|_1$, $\|\cdot\|_2$, $\|\cdot\|_\infty$ and the Frobenius norm, by direct calculation. More generally,

LEMMA 14. *Let* $\|\cdot\|$ *be a matrix norm. Then:*

(i) *if* $J_n^* = \|\mathbf{J}\|^*$, *then* $J_n^* = J_n$;
(ii) $n \le J_n \le C_n$;
(iii) *if* $\|\cdot\|$ *is an operator norm, then* $J_n = C_n$;
(iv) *if* $\|\cdot\|$ *is induced by a vector norm which is symmetric in the coordinates, then* $J_n = n$;
(v) *if* $\|\cdot\|_p$ *is induced by the vector p-norm* $(1 \le p \le \infty)$, *then* $J_n = n$;
(vi) *if* $\|\cdot\|_w = \|W \cdot W^{-1}\|_1$, *where* $W = \operatorname{diag}(w)$ *for a column vector* $w > \mathbf{0}$ *with* $\|w\|_1 = 1$, *then* $J_n = 1/w_{\min}$, *where* $w_{\min} = \min_i w_i$.

PROOF. We have:

(i) $J_n^* = \|\mathbf{J}\|^* = \|\mathbf{J}^{\mathrm{T}}\| = \|\mathbf{J}\| = J_n$.

(ii) $\mathbf{J} = \mathbf{1}\mathbf{1}^{\mathrm{T}}$ so $n\mathbf{1} = \mathbf{J}\mathbf{1}$. Thus $n\|\mathbf{1}\| \le \|\mathbf{J}\|\|\mathbf{1}\|$. Now $\|\mathbf{1}\| \ne 0$, so cancellation gives the first inequality. The second follows by submultiplicativity and duality.

(iii) $J_n = \|\mathbf{J}\| = \|\mathbf{1}\mathbf{1}^{\mathrm{T}}\| = \max_{x \ne 0} \|\mathbf{1}\mathbf{1}^{\mathrm{T}}x\|/\|x\|$, where $x$ is a length-$n$ vector. Pulling scalar multiples out of the vector norm in the numerator, this is equal to $\|\mathbf{1}\| \max_{x \ne 0} |\mathbf{1}^{\mathrm{T}}x|/\|x\|$. Now by Lemma 5, $\|\mathbf{1}\| = \|1\|\,\|\mathbf{1}\|$, and hence $J_n = \|\mathbf{1}\| \max_{x \ne 0} \|\mathbf{1}^{\mathrm{T}}x\|/\|x\| = \|\mathbf{1}\|\|\mathbf{1}^{\mathrm{T}}\| = C_n$.

(iv) Let $x$ be any column vector such that $\mathbf{1}^{\mathrm{T}}x = n$. Let $x_\sigma$ be $x$ after a coordinate permutation $\sigma$, and $\bar{x} = \sum_\sigma x_\sigma/n!$. Clearly, $\bar{x} = \mathbf{1}$. Also, $\|\bar{x}\| \le \|x\|$, and $\mathbf{1}^{\mathrm{T}}\bar{x} = \mathbf{1}^{\mathrm{T}}x = n$ by subadditivity of $\|\cdot\|$ and symmetry, so $\|\mathbf{1}\|^* = \max_{x \ne 0} \mathbf{1}^{\mathrm{T}}x/\|x\| \le \max_{x \ne 0} \mathbf{1}^{\mathrm{T}}\bar{x}/\|\bar{x}\| = n/\|\mathbf{1}\|$.

(v) This follows directly from (iv).

(vi) $J_n = \|\mathbf{J}\|_w = \|Z\|_1$, where $Z_{ij} = w_i/w_j$. Thus, $J_n = \sum_{i=1}^n w_i / \min_i w_i = 1/w_{\min}$. □

REMARK 1. For an arbitrary matrix norm, we can have $C_n > J_n$. This is true even if the norm is invariant under row and column permutations. For example, $\|\cdot\| = \max\{\|\cdot\|_1, \|\cdot\|_\infty\}$ is a matrix norm, with $\|\mathbf{J}\| = \|\mathbf{1}\| = \|\mathbf{1}\|^* = n$, which even satisfies $\|I\| = 1$ (see [23], page 308). For this norm, $C_n/J_n = n$. In general, the ratio is unbounded, even for a fixed $n$. Consider, for example, $\|\cdot\| = \max\{\|W \cdot W^{-1}\|_1, \|W \cdot W^{-1}\|_\infty\}$, where $W = \mathrm{diag}(v)$ for a column vector $v > \mathbf{0}$ with $\|v\|_1 = 1$. It is easy to show that this is a matrix norm with $C_n/J_n = \max_i v_i / \min_i v_i$, which can be arbitrarily large.

We will use the following technical lemma, which appears as Lemma 9 in [14] for the norm $\|\cdot\|_1$. We show that, for any nonnegative matrix $R$ with $\|R\| < 1$, there is a row vector $w$ which approximately satisfies the condition of Lemma 7, and has $w_{\min}$ not too small.

LEMMA 15. *Let $R \in \mathbb{M}_n^+$, and let $\|\cdot\|$ be a matrix norm such that $\|R\| \leq \mu < 1$. Then for any $0 < \eta < 1 - \mu$, there is a matrix $R' \geq R$ and a row vector $w > \mathbf{0}$ such that $wR' \leq \mu' w$, $\|w\|_\infty = 1$ and $w_{\min} = \min_i w_i \geq \eta/J_n$, where $\mu' = \mu + \eta < 1$.*

PROOF. Let $\mathbf{J}' = \mathbf{J}/J_n$, and $R' = R + \eta\mathbf{J}'$. Then $R'$ is irreducible, and $\|R'\| \leq \|R\| + \eta$. Then by Lemma 10, $\lambda(R') \leq \mu + \eta = \mu'$. Thus, by Lemma 7, there exists $w > \mathbf{0}$ such that $wR' \leq \mu' w$. We normalize so that $\|w\|_\infty = 1$. Then $w \geq \mu' w \geq wR' \geq \eta w\mathbf{J}' = \eta \mathbf{1}^{\mathrm{T}}/J_n$, and hence $w_{\min} \geq \eta/J_n$. □

2.2. *Random update and systematic scan Glauber dynamics.* The framework and notation is from [14, 15]. The set of sites of the spin system will be $V = [n] = \{1, 2, \ldots, n\}$, and the set of spins will be $\Sigma = [q]$. A *configuration* (or *state*) is an assignment of a spin to each site, and $\Omega^+ = \Sigma^n$ denotes the set of all such configurations. Let $M = q^n = |\Sigma|^n = |\Omega^+|$, and we will suppose $\Omega^+ = [M]$.

Local interaction between sites specifies the relative likelihood of possible (local) subconfigurations. Taken together, these give a well-defined probability distribution $\pi$ on the set of configurations $\Omega^+$. *Glauber dynamics* is a Markov chain $(x_t)$ on configurations that updates spins one site at a time, and converges to $\pi$. We measure the convergence of this chain by the *total variation distance* $\mathrm{d_{TV}}(\cdot,\cdot)$. We will abuse notation to write, for example, $\mathrm{d_{TV}}(x_t, \pi)$ rather than $\mathrm{d_{TV}}(\mathcal{L}(x_t), \pi)$. The mixing time $\tau(\varepsilon)$ is then defined by $\tau(\varepsilon) = \min_t\{\mathrm{d_{TV}}(x_t, \pi) \leq \varepsilon\}$. In our setting, $n$ measures the size of configurations in $\Omega^+$, and we presume it to be large. Thus, for convenience, we also use *asymptotic* bounds $\hat{\tau}(\varepsilon)$, which have the property that $\limsup_{n\to\infty} \tau(\varepsilon)/\hat{\tau}(\varepsilon) \leq 1$.

We use the following notation. If $x$ is a configuration and $j$ is a site then $x_j$ denotes the spin at site $j$ in $x$. For each site $j$, $S_j$ denotes the set of pairs of configurations that agree off of site $j$. That is, $S_j$ is the set of pairs $(x, y) \in \Omega^+ \times \Omega^+$ such that, for all $i \neq j$, $x_i = y_i$. For any state $x$ and spin $c$, we use $x \to^j c$ for the state $y$ such that $y_i = x_i$ $(i \neq j)$ and $y_j = c$. For each site $j$, we have a transition matrix $P^{[j]}$ on the state space $\Omega^+$ which satisfies two properties:

(i) $P^{[j]}$ changes one configuration to another by updating only the spin at site $j$. That is, if $P^{[j]}(x, y) > 0$, then $(x, y) \in S_j$.

(ii) The equilibrium distribution $\pi$ is invariant with respect to $P^{[j]}$. That is, $\pi P^{[j]} = \pi$.

*Random update Glauber dynamics* corresponds to a Markov chain $\mathcal{M}^\dagger$ with state space $\Omega^+$ and transition matrix $P^\dagger = (1/n)\sum_{j=1}^n P^{[j]}$. *Systematic scan Glauber dynamics* corresponds to a Markov chain $\mathcal{M}_\rightarrow$ with state space $\Omega^+$ and transition matrix $P_\rightarrow = \prod_{j=1}^n P^{[j]}$.

It is well known that the mixing times $\tau_{\mathrm{r}}(\varepsilon)$ of $\mathcal{M}^\dagger$ and $\tau_{\mathrm{s}}(\varepsilon)$ of $\mathcal{M}_\rightarrow$ can be bounded in terms of the influences of sites on each other. To be more precise, let $\mu_j(x, \cdot)$ be the distribution on spins at site $j$ induced by $P^{[j]}(x, \cdot)$. Thus, $\mu_j(x, c) = P^{[j]}(x, x \to^j c)$. Now let $\hat{\varrho}_{ij}$ be the *influence* of site $i$ on site $j$, which is given by $\hat{\varrho}_{ij} = \max_{(x,y)\in S_i} \mathrm{d}_{\mathrm{TV}}(\mu_j(x, \cdot), \mu_j(y, \cdot))$. A *dependency matrix* for the spin system is any $n \times n$ matrix $R = (\varrho_{ij})$ such that $\varrho_{ij} \geq \hat{\varrho}_{ij}$. Clearly, we may assume $\varrho_{ij} \leq 1$.

Given a dependency matrix $R$, let $\varrho_j$ denote the $j$th column of $R$, for $j \in [n]$. Now let $R_j \in \mathbb{M}_n^+$ be the matrix which is an identity except for column $j$, which is $\varrho_j$, that is,

$$(R_j)_{ik} = \begin{cases} 1, & \text{if } i = k \neq j; \\ \varrho_{ij}, & \text{if } k = j; \\ 0, & \text{otherwise.} \end{cases} \tag{1}$$

Let $R^\dagger = \frac{1}{n}\sum_{j=1}^n R_j = \frac{n-1}{n} I + \frac{1}{n} R$ define the *random update matrix* for $R$, and let $\vec{R} = R_1 R_2 \cdots R_n$ define the *scan update matrix* for $R$.

2.3. *The applicability of Lemmas 1 and 2.* In this section, we give an example of a family of spin systems for which Lemmas 1 and 2 can be used to prove rapid mixing, while previous theorems are inapplicable.

*Facilitated spin models* (see [5]) are a class of spin systems in which each spin is either resampled from its equilibrium distribution or is not resampled, depending on whether the surrounding configuration satisfies a local constraint. Consider the following variant of a facilitated spin model on $n$ sites. On each step of the dynamics, a site $j$ is chosen uniformly at random. The spin at the site is sampled from its equilibrium distribution, which is the uniform distribution on $\{0, 1\}$, except that, if any of sites $j - 2$, $j - 1$, or $j+1$ has spin 1, then the resampling only occurs with probability $\delta$ for some $\delta \in (0,1)$.

Let $M$ be the $n \times n$ matrix which has a 1 in entries $(i,i-1)$, $(i,i+1)$, and $(i,i+2)$ (for $i \in \{1,\dots,n\}$) and 0 in all other entries. The dependency matrix $R$ of the spin system is $\frac{1-\delta}{2}M$.

Now suppose, for example, that for $n > 15$, we choose $\delta = \frac{1}{3} - \frac{4}{3(n-2)}$. For these parameters, we will show that the $L_1$, $L_\infty$ and spectral norms of $R$ are all at least 1 (and the $L_\infty$ norm exeeds 1) but $\lambda(R) < 1$. We can draw the following conclusions.

- Since the $L_1$ norm of $R$ is at least 1, the Dobrushin-condition methods of [11, 19, 35, 39] cannot be used to show that random update Glauber dynamics or systematic scan Glauber dynamics are rapidly mixing for this spin system.
- Since the $L_\infty$ norm of $R$ exceeds 1, the methods of [4, 14] cannot be used to show that random update Glauber dynamics or systematic scan Glauber dynamics is rapidly mixing.
- Since the spectral norm of $R$ is at least 1, the methods of [22] are not applicable.
- However, since $\lambda(R) < 1$, by Lemma 11, there is a norm $\|\cdot\|$ with $\|R\| < 1$ and Lemmas 1 and 2 can be used to show rapid mixing of both random update Glauber dynamics and systematic scan Glauber dynamics.

Here is a proof that the $L_1$, $L_\infty$ and spectral norms of $R$ are all at least 1 (and the $L_\infty$ norm exceeds 1) but $\lambda(R) < 1$ (as claimed above).

Let $b = 2/(1-\delta) = 3(1-2/n)$. Each norm of $R$ is the corresponding norm of $M$ divided by $b$, so we wish to show that the $L_1$, $L_\infty$, and spectral norms of $M$ are at least $b$, but that $\lambda(M) < b$.

The $L_1$ and $L_\infty$ norms are easy, so start with the spectral norm $\|M\|_2 = \sqrt{\lambda(P)}$, where $P$ denotes $M^T M$. Since $P$ is symmetric, by Lemma 12, $\lambda(P) = \nu(P)$. Let $x$ be the length-$n$ vector in which every entry is $1/\sqrt{n}$. Then $\nu(P) \geq x^T P x = (1/n)\sum_{i,j} P_{i,j} = (1/n)(9n-18) = 9 - 18/n$. Thus, $\|M\|_2 \geq 3(1-2/n)^{1/2} \geq b$.

Finally, we wish to show $\lambda(M) < b$. By Lemma 7, it suffices to find $w > \mathbf{0}$ satisfying $Mw \leq \mu w$. This will imply $\lambda(M) \leq \mu$. We will take $\mu = 2.62$ which is less than $b$ for $n > 15$. Let $x = 1.525$ and $w_j = x^{-j}$. Then the $i$th row of $Mw$ is at most

$$\begin{aligned} w_{i-1} + w_{i+1} + w_{i+2} &= x^{-i+1} + x^{-i-1} + x^{-i-2} \\ &= x^{-i}\left(x + \frac{1}{x} + \frac{1}{x^2}\right) < w_i \mu, \end{aligned}$$

so we are finished.

**3. Mixing conditions for Glauber dynamics.** There are two approaches to proving mixing results based on the dependency matrix, *path coupling* and *Dobrushin uniqueness*. These are, in a certain sense, dual to each other. All the results given here can be derived equally well using either approach, as we will show.

3.1. *Path coupling.* First, consider applying path coupling to the random update Glauber dynamics. We will begin by proving a simple property of $R^\dagger$.

LEMMA 16. *Let $R$ be a dependency matrix for a spin system, and $\|\cdot\|$ any operator norm such that $\|R\| \le \mu < 1$. Then $\|R^\dagger\| \le \mu^\dagger$ where $\mu^\dagger = 1 - \frac{1}{n}(1-\mu) < 1$.*

PROOF. $\|R^\dagger\| \le \frac{n-1}{n}\|I\| + \frac{1}{n}\|R\| = \frac{n-1}{n} + \frac{1}{n}\|R\| \le 1-(1-\mu)/n = \mu^\dagger$. □

We can use this to bound the mixing time of the random update Glauber dynamics.

LEMMA 17. *Suppose $R$ is a dependency matrix for a spin system, and let $\|\cdot\|$ be any matrix norm. If $\|R\| \le \mu < 1$, then the mixing time $\tau_r(\varepsilon)$ of random update Glauber dynamics is at most $n(1-\mu)^{-1}\ln(C_n/\varepsilon)$.*

PROOF. We will use path coupling. See, for example, [18]. Let $x_0, y_0 \in \Omega^+$ be the initial configurations of the coupled chains, and $x_t, y_t$ the states after $t$ steps of the coupling. The *path* $z_0, \dots, z_n$ from $x_t$ to $y_t$ has states $z_0 = x_t$, and $z_i = (z_{i-1} \to^i y_t(i))$ $(i \in [n])$, so $z_n = y_t$.

We define a distance metric between pairs of configurations as follows. For every $i \in [n]$, we choose a constant $0 < \delta_i \le 1$, and we define the distance between configurations in $S_i$ to be $\delta_i$. That is, for every $(x,y) \in S_i$, we define $\mathrm{d}_\delta(x,y) = \delta_i$. We then lift these distances to a path metric. In particular, for every pair of configurations $(x,y)$, $\mathrm{d}_\delta(x,y) = \sum_{i=1}^n \delta_i I\{x(i) \ne y(i)\}$. The $\delta_i$ $(i \in [n])$ make up a column vector $\delta > \mathbf{0}$. Note that $\mathrm{d}_{\mathbf{1}}(\cdot,\cdot)$ is the usual *Hamming distance.*

Following the path-coupling paradigm, we now define a coupling of one step for each pair of starting states in $S_i$ (for every $i \in [n]$). This gives us a coupling of one step for every pair $(z_i, z_{i+1})$ in the path between $x_t$ and $y_t$ and these can be composed to obtain a coupling of one step from the starting pair $(x_t, y_t)$.

The coupling will be to make the same vertex choice for all $(x_t, y_t) \in S_i$ and then maximally couple the spin choices. With this coupling, $\varrho_{ij}$ bounds the probability of creating a disagreement at site $j$ for any $(x_t, y_t) \in S_i$ and time $t$.

Now consider an arbitrary pair of configurations $(x_t, y_t)$. Let $\beta_t(i) = \Pr(x_t(i) \neq y_t(i))$ determine a row vector $\beta_t$, so $\mathbb{E}[\mathrm{d}_\delta(x_t, y_t)] = \beta_t\delta$. Clearly, $\mathbf{0} \leq \beta_t \leq \mathbf{1}^{\mathrm{T}}$. Since $\beta_t(i)$ and $\Pr(x_{t+1} = x, y_{t+1} = y \,|\, x_t, y_t)$ are independent, it follows that

$$(2)\quad \beta_{t+1}\delta = \mathbb{E}[\mathrm{d}_\delta(x_{t+1}, y_{t+1})] \;\leq\; \sum_{i=1}^{n} \beta_t(i)\left(\delta_i - \frac{\delta_i}{n} + \sum_{j=1}^{n} \frac{\delta_j \varrho_{ij}}{n}\right) = \beta_t R^{\dagger}\delta.$$

[The $i$th term in the sum comes from considering how the distance between $z_{i-1}$ and $z_i$ changes under the coupling. Assuming $z_{i-1}$ and $z_i$ differ (at site $i$) then $\delta_i$ is the reduction in distance that comes about by updating site $i$ and removing the disagreement there, while $\delta_j \varrho_{ij}$ is the expected increase in distance that arises when site $j$ is updated and a disagreement is created there.] Now equation (2) holds for all $\delta$ with $0 < \delta_i \leq 1$. In particular, for any $\varepsilon$, it holds for any vector $\delta$ in which one component is 1 and the other components are $\varepsilon$. Taking the limit, as $\varepsilon \to 0$, we find that componentwise,

$$(3)\qquad \beta_{t+1} \leq \beta_t R^{\dagger}.$$

Now, using (3) and induction on $t$, we find that

$$(4)\qquad \beta_{t+1} \leq \beta_0 {R^{\dagger}}^{t+1}.$$

Equation (4) implies $\beta_{t+1}\mathbf{1} \leq \beta_0 {R^{\dagger}}^{t+1}\mathbf{1}$. Using the coupling lemma [13, 29],

$$\begin{aligned} \mathrm{d}_{\mathrm{TV}}(x_t, y_t) &\leq \Pr(x_t \neq y_t) \leq \sum_{i=1}^{n} \Pr(x_t(i) \neq y_t(i)) \\ &= \beta_t\mathbf{1} \leq \beta_0 {R^{\dagger}}^{t}\mathbf{1} \;\leq\; \mathbf{1}^{\mathrm{T}} {R^{\dagger}}^{t}\mathbf{1}. \end{aligned}$$

Now applying Lemma 4 with $c = \mathbf{1}^{\mathrm{T}} {R^{\dagger}}^{t}\mathbf{1}$ and using submultiplicativity [property (iv) of matrix norms],

$$\mathbf{1}^{\mathrm{T}} {R^{\dagger}}^{t}\mathbf{1} \;\leq\; \|\mathbf{1}\|\|R^{\dagger}\|^{t}\|\mathbf{1}^{\mathrm{T}}\| = C_n\|R^{\dagger}\|^{t}.$$

But $\|R^{\dagger}\| \leq \mu^{\dagger} = 1 - (1-\mu)/n$ by Lemma 16. Thus, when $t \geq n(1-\mu)^{-1}\ln(C_n/\varepsilon)$,

$$\mathrm{d}_{\mathrm{TV}}(x_t, y_t) \;\leq\; C_n {\mu^{\dagger}}^{t} \;=\; C_n(1-(1-\mu)/n)^{t} \;\leq\; C_n e^{-t(1-\mu)/n} \;\leq\; \varepsilon. \qquad \square$$

Corollary 18. *Let $R$ be a dependency matrix for a spin system. Then the mixing time $\tau_r(\varepsilon)$ of random update Glauber dynamics is at most $n(1-\mu)^{-1}\ln(n/\varepsilon)$ if $R$ satisfies any of the following:*

(i) *the Dobrushin condition $\alpha = \|R\|_1 \leq \mu < 1$;*
(ii) *the Dobrushin-Shlosman condition $\alpha' = \|R\|_\infty \leq \mu < 1$;*
(iii) *a $p$-norm condition $\|R\|_p \leq \mu < 1$ for any $1 < p < \infty$.*

PROOF. This follows directly from Lemma 17, since $C_n = J_n = n$ for these norms, by Lemma 14. □

COROLLARY 19. *Let $R$ be a dependency matrix for a spin system. Suppose either of the following conditions holds:*

(i) *$w > \mathbf{0}$ is a row vector such that $wR \le \mu w$, $\|w\|_\infty = 1$ and $w_{\min} = \min_i w_i$;*

(ii) *$w > \mathbf{0}$ is a column vector such that $Rw \le \mu w$, $\|w\|_1 = 1$ and $w_{\min} = \min_i w_i$.*

*Then the mixing time $\tau_r(\varepsilon)$ of random update Glauber dynamics is at most $n(1-\mu)^{-1}\ln(1/w_{\min}\varepsilon)$.*

PROOF. Both are proved similarly, using Lemma 17 with a suitable operator norm, so $C_n = J_n$.

(i) Let $W = \operatorname{diag}(w)$ define the norm $\|R\|_w = \|WRW^{-1}\|_1$. Then $\|R\|_w \le \mu$, and $J_n = 1/w_{\min}$ by Lemma 14.

(ii) Let $W = \operatorname{diag}(w)$ define the norm $\|R\|_w = \|W^{-1}RW\|_\infty = \|WR^{\mathrm{T}}W^{-1}\|_1$. Then $\|R\|_w \le \mu$, and $J_n = 1/w_{\min}$ by Lemma 14. □

In the setting of Corollary 19(i), we can also show contraction of the associated metric $\mathrm{d}_w(\cdot,\cdot)$.

LEMMA 20. *Suppose $R$ is a dependency matrix for a spin system, and let $w > \mathbf{0}$ be a column vector such that $Rw \le \mu w$. Then $\mathbb{E}[\mathrm{d}_w(x_{t+1}, y_{t+1})] \le \mu^\dagger \mathbb{E}[\mathrm{d}_w(x_t, y_t)]$ for all $t \ge 0$.*

PROOF. Note that $R^\dagger w = \frac{n-1}{n}w + \frac{1}{n}Rw \le (\frac{n-1}{n} + \frac{1}{n}\mu)w = \mu^\dagger w$. Putting $\delta = w$ in (2),

$$\mathbb{E}[\mathrm{d}_w(x_{t+1}, y_{t+1})] = \beta_{t+1}w \le \beta_t R^\dagger w \le \mu^\dagger \beta_t w = \mu^\dagger \mathbb{E}[\mathrm{d}_w(x_t, y_t)]. \qquad \square$$

REMARK 2. We may be able to use Lemma 20 obtain a polynomial mixing time in the "equality case" $\mu^\dagger = 1$ of path coupling. However, it is difficult to give a general result other than in "soft core" systems, where all spins can be used to update all sites in every configuration. See [3] for details. We will not pursue this here, however. Note that mixing for the equality case apparently cannot be obtained from the Dobrushin analysis of Section 3.2. This is perhaps the most significant difference between the two approaches.

We would like to use an eigenvector in Corollary 19, since then $\mu = \lambda(R) \le \|R\|$ for any norm. An important observation is that we cannot necessarily do this because $R$ may not be irreducible (so $w_{\min}$ may be 0) or $w_{\min}$ may simply be too small.

$$R = \frac{1}{10}\begin{bmatrix} 0 & 1 & 0 & 0 & 0 & \cdots & 0 & 0 & 0 \\ 8 & 0 & 1 & 0 & 0 & \cdots & 0 & 0 & 0 \\ 0 & 4 & 0 & 1 & 0 & \cdots & 0 & 0 & 0 \\ 0 & 0 & 4 & 0 & 1 & \cdots & 0 & 0 & 0 \\ \vdots & \vdots & & \ddots & \ddots & \ddots & & \vdots & \vdots \\ 0 & 0 & 0 & \cdots & 4 & 0 & 1 & 0 & 0 \\ 0 & 0 & 0 & \cdots & 0 & 4 & 0 & 1 & 0 \\ 0 & 0 & 0 & \cdots & 0 & 0 & 4 & 0 & 2 \\ 0 & 0 & 0 & \cdots & 0 & 0 & 0 & 4 & 0 \end{bmatrix},$$

$$\text{that is,}\ \ \rho_{ij} = \begin{cases} 0.1, & 1 \le i \le n-2,\ j = i+1; \\ 0.4, & 3 \le i \le n,\ j = i-1; \\ 0.8, & i = 2,\ j = 1; \\ 0.2, & i = n-1,\ j = n; \\ 0, & \text{otherwise.} \end{cases}$$

FIG. 1. *Example 1.*

EXAMPLE 1. Consider the matrix of Figure 1. Here $R$ is irreducible, with $\lambda(R) = 0.4$ and left eigenvector $w$ such that $w_i \propto 2^{-i}$ $(i \in [n])$. Thus, $w_{\min} < w_n/w_1 = 2^{1-n}$ is exponentially small, and Corollary 19(i) would give a mixing time estimate of $\Theta(n^2)$ site updates. In fact, $R$ satisfies the Dobrushin condition with $\alpha = 0.8$ and the Dobrushin–Shlosman condition with $\alpha' = 0.9$, so we know mixing actually occurs in $O(n \log n)$ updates.

However, if we know $\|R\| < 1$ for any norm $\|\cdot\|$, we can use Lemma 15 to create a better lower bound on $w_{\min}$. We apply this observation as follows.

COROLLARY 21. *Let $R$ be a dependency matrix for a spin system, and let $\|\cdot\|$ be any matrix norm. Suppose $\|R\| \le \mu < 1$. Then for any $0 < \eta < 1-\mu$, the mixing time of random update Glauber dynamics is bounded by*

$$\tau_{\mathrm{r}}(\varepsilon) \le n(1-\mu-\eta)^{-1} \ln(J_n/\eta\varepsilon).$$

PROOF. Choose $0 < \eta < 1-\mu$. Let $R'$ be the matrix from Lemma 15. Since $R' \ge R$, it is a dependency matrix for the spin system. Let $w$ be the vector from Lemma 15. Now by Corollary 19, the mixing time is bounded by $\tau_{\mathrm{r}}(\varepsilon) \le n(1-\mu')^{-1} \ln(1/w_{\min}\varepsilon)$. where $w_{\min} \ge \eta/J_n$ and $\mu' = \mu + \eta$. $\square$

From this we can now prove Lemma 1, which is a strengthening of Lemma 17 for an arbitrary norm.

LEMMA 1. *Let $R$ be a dependency matrix for a spin system, and let $\|\cdot\|$ be any matrix norm such that $\|R\| \le \mu < 1$. Then the mixing time of random update Glauber dynamics is bounded by*

$$\hat{\tau}_{\mathrm{r}}(\varepsilon) \sim n(1-\mu)^{-1}\ln((1-\mu)^{-1}J_n/\varepsilon).$$

PROOF. Choose $\eta = (1-\mu)/\ln n$. Substituting this into the mixing time from Corollary 21 now implies the conclusion, since $J_n \ge n$. □

REMARK 3. The mixing time estimate is $\hat{\tau}_{\mathrm{r}}(\varepsilon) \sim n(1-\mu)^{-1}\ln((1-\mu)^{-1}J_n/\varepsilon)$. If $(1-\mu)$ is not too small, for example, if $(1-\mu) = \Omega(\log^{-k} n)$ for any constant $k \ge 0$, we have $\hat{\tau}_{\mathrm{r}}(\varepsilon) \sim n(1-\mu)^{-1}\ln(J_n/\varepsilon)$. Thus, we lose little asymptotically using Lemma 1, which holds for an arbitrary matrix norm, from the mixing time estimate $\hat{\tau}_{\mathrm{r}}(\varepsilon) = n(1-\mu)^{-1}\ln(J_n/\varepsilon)$, which results from applying Corollary 17 with an operator norm $\|\cdot\|$. The condition $(1-\mu) = \Omega(\log^{-k} n)$ holds, for example, when $(1-\mu)$ is a small positive constant, which is the case in many applications.

We can easily extend the analysis above to deal with systematic scan. Here the mixing time $\tau_{\mathrm{s}}(\varepsilon)$ will be bounded as a number of complete scans. The number of individual Glauber steps is then $n$ times this quantity. The following lemma modifies the proof technique of [14], Section 7.

LEMMA 22. *Let $R$ be a dependency matrix for a spin system, and $\|\cdot\|$ any matrix norm. If $\|\vec{R}\| \le \mu < 1$, the mixing time $\tau_{\mathrm{s}}(\varepsilon)$ of systematic scan Glauber dynamics is at most $(1-\mu)^{-1}\ln(C_n/\varepsilon)$. If $\|\cdot\|$ is an operator norm, the mixing time is at most $(1-\mu)^{-1}\ln(J_n/\varepsilon)$.*

PROOF. We use the same notation and proof method as in Lemma 17. Consider an application of $P^{[j]}$, with associated matrix $R_j$, as defined in (1). Then it follows that

$$\mathbb{E}[\mathrm{d}(x_1, y_1)] \;\le\; \sum_{i=1}^{n} \beta_0(i)(\delta_i + \delta_j \varrho_{ij}) \;=\; \beta_0 R_j \delta$$

If as before, $\delta_i = 1$ and $\delta_j \to 0$ for $j \ne i$, we have $\Pr(x_1(i) \ne y_1(i)) \le \beta_0 R_j$. Now it follows that $\mathbb{E}[\mathrm{d}(x_n, y_n)] \le \beta_0(\prod_{i=1}^{n} R_j)\delta = \beta_0 \vec{R}\delta$ and $\mathbb{E}[\mathrm{d}(x_{nt}, y_{nt})] \le \beta_0 \vec{R}^t \delta$. Thus, $\Pr(x_{nt}(i) \ne y_{nt}(i)) \le \beta_0 \vec{R}^t(i)$. Hence,

$$\begin{aligned}\mathrm{d}_{\mathrm{TV}}(x_{nt}, y_{nt}) \le \Pr(x_{nt} \ne y_{nt}) &\le \sum_{i=1}^{n} \Pr(x_{nt}(i) \ne y_{nt}(i)) \\ &\le \beta_0 \vec{R}^t \mathbf{1} \le \mathbf{1}^{\mathrm{T}} \vec{R}^t \mathbf{1} \;\le\; \|\vec{R}\|^t \|\mathbf{1}^{\mathrm{T}}\| \|\mathbf{1}\|.\end{aligned}$$

The remainder of the proof is now similar to Lemma 17. □

The following lemma was proved in a slightly different form in [14, Lemma 11]. It establishes the key relationship between $\vec{R}$ and $R$.

LEMMA 23. *Let $R$ be a dependency matrix for a spin system. Suppose $w > \mathbf{0}$ is a row vector, such that $wR \le \mu w$ for some $\mu \le 1$. Then $w\vec{R} \le \mu w$.*

PROOF. Note that for any row vector $z$, $zR_i = [z_1 \cdots z_{i-1} z\varrho_i z_{i+1} \cdots z_n]$. Since $wR \le \mu w \le w$, $w\varrho_i \le w_i$. Now we can show by induction on $i$ that $wR_1 \cdots R_i \le [w\varrho_1 \cdots w\varrho_i w_{i+1} \cdots w_n]$. For the inductive step, $wR_1 \cdots R_i \le zR_i = [z_1 \cdots z_{i-1} z\varrho_i z_{i+1} \cdots z_n]$ where $z = [w\varrho_1 \cdots w\varrho_{i-1} w_i \cdots w_n]$. But then $z \le w$, so $z\varrho_i \le w\varrho_i$ so $zR_i \le [w\varrho_1 \cdots w\varrho_i w_{i+1} \cdots w_n]$. Taking $i = n$, we have $w\vec{R} \le [w\varrho_1 \cdots w\varrho_n] = wR \le \mu w$. $\square$

COROLLARY 24. $\lambda(\vec{R}) \le \lambda(R)$ *and if* $\|R\|_1 \le 1$, $\|\vec{R}\|_1 \le \|R\|_1$.

PROOF. The first statement follows directly from Lemmas 7 and 23. For the second, note that $\mathbf{1}^{\mathrm{T}} R \le \|R\|_1 \mathbf{1}^{\mathrm{T}}$, so $\mathbf{1}^{\mathrm{T}} \vec{R} \le \|R\|_1 \mathbf{1}^{\mathrm{T}}$ by Lemma 23. But this implies $\|\vec{R}\|_1 \le \|R\|_1$. $\square$

We can now apply this to the mixing of systematic scan. First we show, as in [35], that the Dobrushin criterion implies rapid mixing.

COROLLARY 25. *Let $R$ be a dependency matrix for a spin system. Then if $R$ satisfies the Dobrushin condition $\alpha = \|R\|_1 \le \mu < 1$, the mixing time of systematic scan Glauber dynamics is at most $(1-\mu)^{-1} \ln(n/\varepsilon)$.*

PROOF. This follows from Lemma 22 and Corollary 24, since $J_n = n$ for the norm $\|\cdot\|_1$. $\square$

Next we show, as in [14], Section 3.3, that a weighted Dobrushin criterion implies rapid mixing.

COROLLARY 26. *Let $R$ be a dependency matrix for a spin system. Suppose $w > \mathbf{0}$ is a row vector satisfying $\|w\|_\infty = 1$ and $wR \le \mu w$ for some $\mu < 1$. Let $w_{\min} = \min_i w_i$. Then the mixing time $\tau_S(\varepsilon)$ of systematic scan Glauber dynamics is bounded by $(1-\mu)^{-1} \ln(1/w_{\min}\varepsilon)$.*

PROOF. By Lemma 23, $w\vec{R} \le \mu w$. We use the norm $\|\cdot\|_w = \|W \cdot W^{-1}\|_1$, where $W = \mathrm{diag}(w)$. Then apply Lemma 22 with $\|\vec{R}\|_w \le \mu$. $\square$

Once again, we cannot necessarily apply Corollary 26 directly since $w_{\min}$ may be too small (or even 0). Applying Corollary 26 to Example 1 would give a mixing time estimate of $\Theta(n)$ scans. However, $R$ satisfies the Dobrushin condition with $\alpha = 0.8$ so we know mixing actually occurs in $O(\log n)$ scans. Once again, our solution is to perturb $R$ using Lemma 15.

COROLLARY 27. *Let $R$ be a dependency matrix for a spin system and $\|\cdot\|$ a matrix norm such that $\|R\| \le \mu < 1$. Then for any $0 < \eta < 1 - \mu$, the mixing time of systematic scan Glauber dynamics is bounded by*

$$\tau_S(\varepsilon) \le (1 - \mu - \eta)^{-1} \ln(J_n/\eta\varepsilon).$$

PROOF. Let $R'$ be the matrix and $w$ the vector from Lemma 15. Since $R' \ge R$, it is a dependency matrix for the spin system. Now by Corollary 26, the mixing time satisfies $\tau_{\mathrm{S}}(\varepsilon) \le (1-\mu')^{-1} \ln(1/w_{\min}\varepsilon)$, where $w_{\min} = \min_i w_i \ge \eta/J_n$ and $\mu' = \mu + \eta$. □

We can now use this to prove Lemma 2.

LEMMA 2. *Let $R$ be a dependency matrix for a spin system and $\|\cdot\|$ any matrix norm such that $\|R\| \le \mu < 1$. Then the mixing time of systematic scan Glauber dynamics is bounded by*

$$\hat{\tau}_S(\varepsilon) \sim (1 - \mu)^{-1} \ln((1 - \mu)^{-1} J_n/\varepsilon).$$

PROOF. This follows from Corollary 27 exactly as Lemma 1 follows from Corollary 21. □

REMARK 4. If, for example, $\|\cdot\| = \|\cdot\|_p$, for any $1 < p \le \infty$, $J_n = n$, and we obtain a mixing time $\hat{\tau}_{\mathrm{S}}(\varepsilon) \sim (1 - \mu)^{-1} \ln((1 - \mu)^{-1} n/\varepsilon)$. If in addition, $(1 - \mu) = \Omega(\log^{-k} n)$ for any $k \ge 0$ (as in Remark 3), we have $\hat{\tau}_{\mathrm{S}}(\varepsilon) \sim (1 - \mu)^{-1} \ln(n/\varepsilon)$, which matches the bound from Corollary 25 for the norm $\|\cdot\|_1$. Note that there is a difference from the random update case, since here we do not have a result like Lemma 17 which we can apply directly with any operator norm.

3.2. *Dobrushin uniqueness.* The natural view of path coupling in this setting corresponds to multiplying $R^\dagger$ on the left by a row vector $\beta$, as in Lemma 17. The Dobrushin uniqueness approach corresponds to multiplying $R$ on the right by a column vector $\delta$. As we showed in [14], Section 7, these two approaches are essentially equivalent. However, for historical reasons, the Dobrushin uniqueness approach is frequently used in the statistical physics literature. See, for example, [33, 35]. Therefore, for completeness, we will now describe the Dobrushin uniqueness framework, using the notation of [14].

Recall that $\Omega^+ = [M]$. For any column vector $f \in \mathbb{R}^M$, let $\delta_i(f) = \max_{(x,y)\in S_i} |f(x) - f(y)|$. Let $\delta(f)$ be the column vector given by $\delta(f) = (\delta_1(f), \delta_2(f), \dots, \delta_n(f))$. Thus, $\delta : \mathbb{R}^M \to \mathbb{R}^n$. The following lemma gives the key property of this function.

LEMMA 28 ([14], Lemma 10). *The function $\delta$ satisfies $\delta(P^{[j]}f) \le R_j \delta(f)$.*

PROOF. Suppose $(x,y) \in S_i$ maximizes $|P^{[j]}f(x) - P^{[j]}f(y)|$. Then

$$\begin{aligned}
\delta_i(P^{[j]}f) &= |P^{[j]}f(x) - P^{[j]}f(y)| \\
&= \left| \sum_c f(x \to^j c) P^{[j]}(x, x \to^j c) - \sum_c f(y \to^j c) P^{[j]}(y, y \to^j c) \right| \\
&= \left| \sum_c f(x \to^j c)\mu_j(x,c) - \sum_c f(y \to^j c)\mu_j(y,c) \right| \\
&= \left| \sum_c (f(x \to^j c) - f(y \to^j c))\mu_j(x,c) \right. \\
&\qquad \left. + \sum_c f(y \to^j c)\,(\mu_j(x,c) - \mu_j(y,c)) \right| \\
&\le \sum_c |f(x \to^j c) - f(y \to^j c)|\mu_j(x,c) \\
&\qquad + \left| \sum_c f(y \to^j c)\,(\mu_j(x,c) - \mu_j(y,c)) \right|.
\end{aligned}$$

We will bound the two terms in the last expression separately. First,

$$\begin{aligned}
&\sum_c |f(x \to^j c) - f(y \to^j c)|\mu_j(x,c) \\
&\qquad \le \max_c |f(x \to^j c) - f(y \to^j c)| \le \mathbb{1}_{i \ne j}\delta_i(f).
\end{aligned} \tag{5}$$

For the second, let $f^+ = \max_c f(y \to^j c)$, $f^- = \min_c f(y \to^j c)$ and $f^0 = \frac{1}{2}(f^+ + f^-)$. Note that $f^+ - f^0 = \frac{1}{2}(f^+ - f^-) \le \frac{1}{2}\delta_j(f)$. Then since $\sum_c(\mu_j(x,c) - \mu_j(y,c)) = 0$,

$$\begin{aligned}
&\left| \sum_c f(y \to^j c)\,(\mu_j(x,c) - \mu_j(y,c)) \right| \\
&\qquad = \left| \sum_c \left( f(y \to^j c) - f^0 \right) (\mu_j(x,c) - \mu_j(y,c)) \right|
\end{aligned}$$

$$
\begin{aligned}
(6) \qquad & \le 2\mathrm{d}_{\mathrm{TV}}(\mu_j(x,\cdot),\mu_j(y,\cdot)) \max_c |f(y \to^j c) - f^0| \\
&= 2\mathrm{d}_{\mathrm{TV}}(\mu_j(x,\cdot),\mu_j(y,\cdot))(f^+ - f^0) \\
&\le \varrho_{ij}\delta_j(f).
\end{aligned}
$$

The conclusion now follows by adding (5) and (6). □

The following lemma allows us to apply Lemma 28 to bound mixing times.

LEMMA 29. *Let $\mathcal{M} = (X_t)$ be a Markov chain with transition matrix $P$, and $\|\cdot\|$ a matrix norm. Suppose there is a matrix $R$ such that, for any column vector $f \in \mathbb{R}^M$, $\delta(Pf) \le R\delta(f)$, and $\|R\| \le \mu < 1$. Then the mixing time of $\mathcal{M}$ is bounded by*

$$
\tau(\varepsilon) \le (1-\mu)^{-1} \ln(C_n/\varepsilon)
$$

PROOF. For a column vector $f_0$, let $f_t$ be the column vector $f_t = P^t f_0$. Let $\pi$ be the row vector corresponding to the stationary distribution of $\mathcal{M}$. Note that $\pi f_t = \pi P^t f_0$, which is $\pi f_0$ since $\pi$ is a left eigenvector of $P$ with eigenvalue 1.

Now let $f_0$ be the indicator vector for an arbitrary subset $A$ of $[M] = \Omega^+$. That is, let $f_0(z) = 1$ if $z \in A$ and $f_0(z) = 0$ otherwise. Then since $P^t(x,y) = \Pr(X_t = y \mid X_0 = x)$, we have $f_t(x) = \Pr(X_t \in A \mid x_0 = x)$. Also, $\pi f_t = \pi f_0 = \pi(A)$ for all $t$. Let $f_t^- = \min_z f_t(z)$ and $f_t^+ = \max_z f_t(z)$. Since $\pi$ is a probability distribution, $f_t^- \le \pi f_t \le f_t^+$, so $f_t^- \le \pi(A) \le f_t^+$.

By induction on $t$, using the condition in the statement of the lemma, we have $\delta(f_t) \le R^t\delta(f_0)$. But $R^t\delta(f_0) \le R^t\mathbf{1}$. Now, consider states $x, y$ such that $f_t(x) = f_t^-$, $f_t(y) = f_t^+$. Let $z_i$ $(i = 0, 1, \dots, n)$ be the path of states from $x$ to $y$ used in the proof of Lemma 17. Then

$$
\begin{aligned}
f_t^+ - f_t^- = f_t(y) - f_t(x) &\le \sum_{i=1}^{n} |f_t(z_i) - f_t(z_{i-1})| \\
&\le \sum_{i=1}^{n} \delta_i(f_t) = \mathbf{1}^{\mathrm{T}}\delta(f_t) \le \mathbf{1}^{\mathrm{T}} R^t \mathbf{1}.
\end{aligned}
$$

This implies that $\max_x |\Pr(x_t \in A \mid x_0 = x) - \pi(A)| \le \mathbf{1}^{\mathrm{T}} R^t \mathbf{1}$. Since $A$ is arbitrary, for all $t \ge (1-\mu)^{-1}\ln(C_n/\varepsilon)$ we have

$$
\begin{aligned}
\mathrm{d}_{\mathrm{TV}}(x_t, \pi) &\le \mathbf{1}^{\mathrm{T}} R^t \mathbf{1} \le \|R\|^t \|\mathbf{1}\| \|\mathbf{1}\|^* \\
&= C_n \|R\|^t \le C_n \mu^t \le C_n e^{-(1-\mu)t} \le \varepsilon.
\end{aligned}
$$
□

The following lemma and Lemma 17, whose proof follows, enable us to use Lemma 29 to bound the mixing time of random update Glauber dynamics.

LEMMA 30. *Let $R$ be a dependency matrix for a spin system Let $R^\dagger$ be the random update matrix for $R$. Then for $f \in \mathbb{R}^M$, $\delta(P^\dagger f) \le R^\dagger \delta(f)$.*

PROOF. For each $i \in [n]$, from the definition of $\delta_i$, $\delta_i(f) \ge 0$ and, for any $c \in \mathbb{R}$ and $f \in \mathbb{R}^M$, $\delta_i(cf) = |c|\delta_i(f)$. Also, $\delta_i(f_1 + f_2) \le \delta_i(f_1) + \delta_i(f_2)$ for any $f_1, f_2 \in \mathbb{R}^M$. Now,

$$\delta(P^\dagger f) = \delta\left(\frac{1}{n}\sum_{j=1}^n P^{[j]} f\right) = \frac{1}{n}\delta\left(\sum_{j=1}^n P^{[j]} f\right) \le \frac{1}{n}\sum_{j=1}^n \delta(P^{[j]} f).$$

By Lemma 28, this is at most $\frac{1}{n}\sum_{j=1}^n R_j \delta(f) = R^\dagger \delta(f)$. □

REMARK 5. The proof shows that $\delta_i(f)$ is a (vector) *seminorm* for all $i \in [n]$. It fails to be a *norm* because $\delta_i(f) = 0$ does not imply $f = 0$. For example, $\delta_i(\mathbf{1}) = 0$ for all $i \in [n]$.

We can now give a proof of Lemma 17 using this approach.

LEMMA 17. *Suppose $R$ is a dependency matrix for a spin system, and let $\|\cdot\|$ be any matrix norm. If $\|R\| \le \mu < 1$, then the mixing time $\tau_r(\varepsilon)$ of random update Glauber dynamics is at most $n(1-\mu)^{-1}\ln(C_n/\varepsilon)$.*

PROOF. By Lemma 16, $\|R^\dagger\| \le \mu^\dagger = 1 - \frac{1}{n}(1-\mu)$ and by Lemma 30, $\delta(P^\dagger f) \le R^\dagger \delta(f)$. Then by Lemma 29, $\tau_r(\varepsilon) \le (1-\mu^\dagger)^{-1}\ln(C_n/\varepsilon) = n(1-\mu)^{-1}\ln(C_n/\varepsilon)$. □

Corollaries 18 and 19 and the rest of that section now follow exactly as before. A similar analysis applies to systematic scan, though it is slightly easier. It relies on the analogue of Lemma 30, which in this case is immediate from Lemma 28.

LEMMA 31. *Let $R$ be a dependency matrix for a spin system. Let $\vec{R}$ be the scan update matrix for $R$. Then for any $f \in \mathbb{R}^M$, $\delta(\vec{P} f) \le \vec{R}\delta(f)$.* □

We can now give a proof of Lemma 22 using this approach.

LEMMA 19. *Let $R$ be a dependency matrix for a spin system, and $\|\cdot\|$ any matrix norm. If $\|\vec{R}\| \le \mu < 1$, the mixing time $\tau_s(\varepsilon)$ of systematic scan Glauber dynamics is at most $(1-\mu)^{-1}\ln(C_n/\varepsilon)$. If $\|\cdot\|$ is an operator norm, the mixing time is at most $(1-\mu)^{-1}\ln(J_n/\varepsilon)$.*

PROOF. By Lemma 31, for any $f \in \mathbb{R}^M$, $\delta(\vec{P} f) \le \vec{R}\delta(f)$. Then by assumption, $\|\vec{R}\| \le \mu < 1$. Now apply Lemma 29. □

The results following Lemma 22 in Section 3.1 can then be obtained identically to the proofs given there.

3.3. *Improved analysis of systematic scan.* We may improve the analysis of Corollary 27 for the case in which the diagonal of $R$ is $\mathbf{0}$, which is the case for the *heat bath* dynamics. For $\sigma \geq 0$, define $R^\sigma$ by

$$R^\sigma_{ij} = \begin{cases} \sigma R_{ij}, & \text{if } 1 \leq i < j \leq n; \\ R_{ij}, & \text{otherwise,} \end{cases}$$

so $R^\sigma$ has its upper triangle scaled by $\sigma$. Let $\varrho^\sigma_j$ denote the $j$th column of $R^\sigma$, for $j \in [n]$. We can now prove the following strengthening of Lemma 23.

LEMMA 32. *If $wR^\sigma \leq \sigma w$, for some $w \geq \mathbf{0}$ and $0 \leq \sigma \leq 1$, then $w\vec{R} \leq wR^\sigma$.*

PROOF. We prove by induction that

$$\begin{aligned} wR_1R_2 \cdots R_i &\leq [w\varrho^\sigma_1 \ \cdots \ w\varrho^\sigma_i \ \ w_{i+1} \ \cdots \ w_n] \\ &\leq [\sigma w_1 \ \cdots \ \sigma w_i \ w_{i+1} \ \cdots \ w_n]. \end{aligned}$$

The second inequality follows by assumption. The hypothesis is clearly true for $i = 0$. For $i > 0$,

$$\begin{aligned} wR_1R_2 \cdots R_{i-1}R_i &\leq [w\varrho^\sigma_1 \ \cdots \ w\varrho^\sigma_{i-1} \ w_i \ w_{i+1} \ \cdots \ w_n]R_i \\ &= [w\varrho^\sigma_1 \ \cdots \ w\varrho^\sigma_{i-1} \ \widetilde{w}\varrho_i \ w_{i+1} \ \cdots \ w_n], \end{aligned}$$

where $\widetilde{w} = [w\varrho^\sigma_1 \ \cdots \ w\varrho^\sigma_{i-1} \ w_i \ \cdots \ w_n] \leq [\sigma w_1 \ \cdots \ \sigma w_{i-1} \ w_i \ \cdots \ w_n]$. It follows that $\widetilde{w}\varrho_i \leq w\varrho^\sigma_i$, continuing the induction. Putting $i = n$ gives the conclusion. $\square$

LEMMA 33. *If $R$ is symmetric and $\lambda = \lambda(R) < 1$ then $\lambda(R^\sigma) \leq \sigma$ if $\sigma = \lambda/(2-\lambda)$.*

PROOF. We have $\lambda = \nu(R) = \|R\|_2$ by Lemma 12. Since $R$ is symmetric with zero diagonal, $x^TR^\sigma x = \frac{1}{2}(1+\sigma)x^TRx$. It follows that $\lambda(R^\sigma) \leq \nu(R^\sigma) = \frac{1}{2}(1+\sigma)\nu(R) = \frac{1}{2}(1+\sigma)\lambda$. Therefore, $\lambda(R^\sigma) \leq \sigma$ if $\lambda \leq 2\sigma/(1+\sigma)$. This holds if $\sigma \geq \lambda/(2-\lambda)$. $\square$

LEMMA 34. *Let $R$ be symmetric with zero diagonal and $\|R\|_2 = \lambda(R) = \lambda < 1$, and $0 < \eta < 1 - \lambda$. Let $\mu = \lambda + \eta < 1$. Then the mixing time of systematic scan is at most*

$$\tau_{\mathrm{s}}(\varepsilon) \leq \frac{2-\mu}{2-2\mu} \ln(n/\eta\varepsilon).$$

PROOF. Let $n' = n - 1$ and $S = R + \eta(\mathbf{J} - I)/n'$. Since $S \geq R$, $S$ is a dependency matrix for the original spin system. Also, $S$ is symmetric and its diagonal is 0. Now

$$\lambda(S) = \|S\|_2 = \|R + \eta(\mathbf{J} - I)/n'\|_2 \leq \|R\|_2 + \eta\|\mathbf{J} - I\|_2/n' \leq \lambda + \eta = \mu.$$

Denote by $\vec{S} = S_1 S_2 \cdots S_n$ the scan matrix. Let $\sigma = \mu/(2 - \mu)$. Now by Lemma 33, we have $\lambda(S^\sigma) \leq \sigma$. Furthermore, $S^\sigma$ is irreducible, so by Lemma 7, there exists a row vector $w > \mathbf{0}$ satisfying $wS^\sigma \leq \sigma w$. We can assume without loss of generality that $w$ is normalised so that $\|w\|_\infty = 1$. Finally, we can conclude from Lemma 32 that $w\vec{S} \leq wS^\sigma$.

Since $w\vec{S} \leq wS^\sigma \leq \sigma w$, we have established that convergence is geometric with ratio $\sigma$, but we need a lower bound on $w_{\min} = w_{\min}$ in order to obtain an upper bound on mixing time via Lemma 29. Now

$$\begin{aligned} \sigma w \geq wS^\sigma &\geq w(\sigma R + \sigma\eta(\mathbf{J} - I)/n') \\ &\geq \sigma\eta w(\mathbf{J} - I)/n' = (\sigma\eta/n')(\mathbf{1} - w). \end{aligned}$$

So $w(1 + \eta/n') \geq (\eta/n')\mathbf{1}$, and $w_{\min} \geq \eta/(n' + \eta) \geq \eta/n$. By Corollary 26, the mixing time satisfies

$$\tau_{\mathrm{S}}(\varepsilon) \leq (1 - \sigma)^{-1} \ln(1/w_{\min}\varepsilon) \leq (1 - \sigma)^{-1} \ln(n/\eta\varepsilon). \qquad \square$$

We can now prove Lemma 3.

LEMMA 3. *Let $R$ be symmetric with zero diagonal and $\|R\|_2 = \lambda(R) = \lambda < 1$. Then the mixing time of systematic scan is at most*

$$\hat{\tau}_S(\varepsilon) \sim (1 - \tfrac{1}{2}\lambda)(1 - \lambda)^{-1} \ln((1 - \lambda)^{-1} n/\varepsilon).$$

PROOF. We apply Lemma 34 with $\eta = (1 - \lambda)/\ln n$, and hence $\mu \sim \lambda$. $\square$

REMARK 6. If, as in Remark 3, if $(1 - \lambda) = \Omega(\log^{-k} n)$ for some $k \geq 0$, then we have mixing time $\hat{\tau}_{\mathrm{s}}(\varepsilon) \sim (1 - \frac{1}{2}\lambda)(1 - \lambda)^{-1} \ln(n/\varepsilon)$ for systematic scan. We may compare the number of Glauber steps $n\tau_{\mathrm{s}}(\varepsilon)$ with the estimate $\hat{\tau}_{\mathrm{r}}(\varepsilon) = (1 - \lambda)^{-1} n \ln(n/\varepsilon)$ for random update Glauber dynamics obtained from Corollary 18 using the minimum norm $\|\cdot\|_2$. The ratio is $(1 - \frac{1}{2}\lambda) < 1$. This is close to $\frac{1}{2}$ when $\lambda(R)$ is close to 1, as in many applications.

EXAMPLE 2. Consider coloring a $\Delta$-regular graph with $(2\Delta + 1)$ colors ([24, 33]) using heat bath Glauber dynamics, we have $\lambda(R) = \Delta/(\Delta + 1)$. (See Section 4). Then $(1 - \lambda) = 1/(\Delta + 1) = \Omega(1)$, if $\Delta = O(1)$, and the above ratio is $(1 - \frac{1}{2}\lambda) = (\Delta + 2)/(2\Delta + 2)$. This is close to $\frac{1}{2}$ for large $\Delta$.

Although the improvement in the mixing time bound is a modest constant factor, this provides some evidence in support of the conjecture that systematic scan mixes faster than random update, for Glauber dynamics at least. The improvement is because we know, later in the scan, that most vertices have already been updated. In a random update, some vertices are updated many times before others are updated at all. Lemma 34 suggests that this may be wasteful.

**4. Coloring sparse graphs.** In this section, we consider an application of the methods developed above to graph coloring problems, particularly in *sparse* graphs. By sparse, we will mean here that the number of edges of the graph is at most linear in its number of vertices.

Let $G = (V, E)$, with $V = [n]$, be an undirected (simple) graph or multigraph, without self-loops. Then $d_v$ will denote degree of vertex $v \in V$. If $S \subseteq V$, we will denote the induced subgraph by $G_S = (S, E_S)$. The (symmetric) *adjacency matrix* $A(G)$ is a nonnegative integer matrix, with zero diagonal, giving the number of edges between each pair of vertices. We write $A$ for $A(G)$ and $\lambda(G)$ for $\lambda(A(G))$. Thus, the adjacency matrix of a graph is a 0–1 matrix. We also consider digraphs and directed multigraphs $\vec{G} = (V, \vec{E})$. We denote the indegree and outdegree of $v \in V$ by $d_v^-, d_v^+$, respectively.

If $G$ is a graph with maximum degree $\Delta$, we consider the heat bath Glauber dynamics for properly coloring $V$ with $q > \Delta$ colors. The dependency matrix $R$ for this chain satisfies $\varrho_{ij} \le 1/(q - d_j)$ $(i, j \in [n])$ (see Section 5.2 of [14]). Thus, $R = AD$, where $D = \operatorname{diag}(1/(q - d_j))$. Let $D^{1/2} = \operatorname{diag}(1/\sqrt{q - d_j})$ and $\hat{A} = D^{1/2} A D^{1/2}$. Note that $\hat{A}$ is symmetric. Also, $\lambda(\hat{A}) = \lambda(AD)$, since $(D^{1/2} A D^{1/2})(D^{1/2} x) = \lambda(D^{1/2} x)$ if and only if $ADx = \lambda x$. If $(i, j) \in E$, we have $\hat{A}_{ij} = 1/\sqrt{(q - d_i)(q - d_j)}$. Since $\hat{A} \le \frac{1}{q-\Delta} A$, we have $\lambda(\hat{A}) \le \frac{1}{q-\Delta} \lambda(A)$ from Lemma 9. So if $q > \Delta + \lambda(A)$, we can use Lemmas 1 and 2 to show that scan and Glauber both mix rapidly. For very nonregular graphs, we may have $\lambda(\hat{A}) \ll \frac{1}{q-\Delta} \lambda(A)$. However, $\lambda(\hat{A})$ seems more difficult to estimate than $\lambda(A)$, since it depends more on the detailed structure of $G$. Therefore, we will use the bound $\frac{1}{q-\Delta} \lambda(A)$ in the remainder of this section, and restrict most of the discussion to $\lambda(G)$. The following is well known.

LEMMA 35. *If $G$ has maximum degree $\Delta$ and average degree $\bar{d}$, then $\bar{d} \le \lambda(G) \le \Delta$. If either bound is attained, there is equality throughout and $G$ is $\Delta$-regular.*

PROOF. The vertex degrees of $G$ are the row or column sums of $A(G)$. The upper bound then follows from $\lambda(G) \le \|A\|_1 = \max_{v \in V} d_v = \Delta$ using Lemma 10. For the lower bound, since $G$ is undirected, $\lambda(G) = \nu(A) \ge$ $\frac{\mathbf{1}^{\mathrm{T}}}{\sqrt{n}} A \frac{\mathbf{1}}{\sqrt{n}} = 2|E|/n = \bar{d}$, using Lemma 12. If the lower bound is attained, then the inequalities in the previous line are equalities, so $\mathbf{1}$ is an eigenvector of $A$. Thus, $A\mathbf{1} = \bar{d}\mathbf{1}$, and every vertex has degree $\bar{d} = \Delta$. When the upper bound is attained, since the columns sums of $A$ are at most $\Delta$, $\mathbf{1}A \leq \Delta\mathbf{1} = \lambda\mathbf{1}$, so $\mathbf{1}$ is an eigenvector from Lemma 7 and $\mathbf{1}A = \Delta\mathbf{1}$. Then every vertex has degree $\Delta = \bar{d}$. □

Thus, the resulting bound for coloring will be $q > 2\Delta$ when $G$ is $\Delta$-regular, as already shown by Jerrum [24] or Salas and Sokal [33]. Thus, we can only achieve mixing for $q \leq 2\Delta$ by this approach if the degree sequence of $G$ is nonregular.

We now derive a bound on $\lambda(R)$ for symmetric $R$ which is very simple, but nonetheless can be used to provide good estimates in some applications.

LEMMA 36. *Suppose $R \in \mathbb{M}_n^+$, and we have $R = B + B^{\mathrm{T}}$, for some $B \in \mathbb{M}_n$. If $\|\cdot\|$ is any matrix norm, then $\lambda(R) \leq 2\sqrt{\|B\|\|B\|^*}$.*

PROOF. $\lambda(R) = \|B + B^{\mathrm{T}}\|_2 \leq \|B\|_2 + \|B^{\mathrm{T}}\|_2 = 2\|B\|_2 \leq 2\sqrt{\|B\|\|B\|^*}$, using the self-duality of $\|\cdot\|_2$ and Lemmas 12 and 13. □

COROLLARY 37. *If $R = B + B^{\mathrm{T}}$, then $\lambda(R) \leq 2\sqrt{\|B\|_1\|B\|_\infty}$.*

We can use Corollary 37 as follows. If $R \in \mathbb{M}_n^+$, let $\kappa(R) = \max_{I \subseteq [n]} \sum_{i,j \in I} \varrho_{ij}/2|I|$. We call $\kappa(R)$ the *maximum density* of $R$. Note that $\kappa(R) \geq \frac{1}{2}\max_{i \in [n]} \varrho_{ii}$. Thus, the maximum density $\kappa(G)$ of $A(G)$ for a graph or multigraph $G = (V, E)$ is $\max_{S \subseteq V} |E_S|/|S|$, according with common usage. This measure will be useful for *sparse* graphs. Note that the maximum density can be computed in polynomial time [21]. Note also that, for symmetric $R \in \mathbb{M}_n^+$, the maximum density is a discrete version of the largest eigenvalue, since

$$\kappa(R) = \max_{x \in \{0,1\}^n} \frac{x^{\mathrm{T}} R x}{x^{\mathrm{T}} x} \leq \max_{x \in \mathbb{R}^n} \frac{x^{\mathrm{T}} R x}{x^{\mathrm{T}} x} = \nu(R) = \lambda(R).$$

Also, $\alpha(R) = \|R\|_1 \geq 2\kappa(R)$, since

$$\begin{aligned} \kappa(R) = \max_{I \subseteq [n]} \sum_{i,j \in I} \varrho_{ij}/2|I| &\leq \max_{I \subseteq [n]} \sum_{i \in [n], j \in I} \varrho_{ij}/2|I| \\ &\leq \max_{I \subseteq [n]} \alpha(R)|I|/2|I| = \alpha(R)/2. \end{aligned}$$

We may easily bound the maximum density for some classes of graphs.

For any $a, b \in \mathbb{Z}$, let us define $\mathcal{G}(a,b)$ to be the maximal class of graphs such that:

(i) $\mathcal{G}(a,b)$ is *hereditary* (closed under taking induced subgraphs);
(ii) for all $G = (V,E) \in \mathcal{G}(a,b)$ with $|V| = n$, we have $|E| \leq an - b$.

LEMMA 38. *Let $G \in \mathcal{G}(a,b)$ with $|V| = n$. If:*

(i) *$b \geq 0$, then $\kappa(G) \leq a - b/n$;*
(ii) *$b \leq 0$, let $k^* = a + \frac{1}{2} + \sqrt{(a+\frac{1}{2})^2 - 2b}$, then $\kappa(G) \leq \kappa^* = \max\{(\lfloor k^* \rfloor - 1)/2, a - b/\lceil k^* \rceil\}$.*

PROOF. In case (i), clearly $|E|/|V| \leq a - b/n$. If $S \subset V$, $|E_S|/|S| \leq a - b/|S| \leq a - b/n$. In case (ii), note that $\kappa(G) \leq \frac{1}{n}\binom{n}{2} = \frac{1}{2}(n-1)$ for any simple graph $G$ on $n$ vertices. Thus,

$$\kappa(G) \;\leq\; \max_{1 \leq |S| \leq n} \min\{(|S|-1)/2, a - b/|S|\}.$$

Note that $(s-1)/2$ is increasing in $s$ and $a - b/s$ is decreasing in $s$. Also, $s = k^*$ is the positive solution to $(s-1)/2 = a - b/s$. The other solution is not positive since $b \leq 0$. Thus

$$\kappa(G) \;\leq\; \max\{(\lfloor k^* \rfloor - 1)/2, a - b/\lceil k^* \rceil\} \;= \kappa^*. \qquad \square$$

REMARK 7. We could consider a more general class $\mathcal{G}(a_n, b_n)$, where $|b_n| = o(na_n)$. This includes, for example, subgraphs of the $d$-dimensional hypercubic grid with vertex set $V = [k]^d$ in which each interior vertex has $2d$ neighbors. Then $|E| \leq dn - dn^{1-1/d}$, so $a_n = d$ and $b_n = dn^{1-1/d}$. However, we will not pursue this further here.

We can apply Lemma 38 directly to some classes of sparse graphs.

For the definition of the *tree-width* $t(G)$ of a graph $G$, see [10]. We say that a graph $G$ *has genus* $g$ if it can be embedded into a surface of genus $g$. See [7] for details, but note that that text (and several others) define the genus *of* the graph to be the *smallest* genus of all surfaces in which $G$ can be embedded. We use our definition because it is appropriate for hereditary classes. Thus, for us a planar graph has genus 0, and a graph which can be embedded in the torus has genus 1 (whether or not it is planar).

LEMMA 39. *If a graph $G = (V,E)$ is:*

(i) *a nonregular connected graph with maximum degree $\Delta$, then $G \in \mathcal{G}(\Delta/2, 1)$;*
(ii) *a forest, then $G \in \mathcal{G}(1,1)$;*
(iii) *a graph of tree-width $t$, then $G \in \mathcal{G}(t, t(t+1)/2)$;*

(iv) *a planar graph, then* $G \in \mathcal{G}(3,6)$*;*
(v) *a graph of genus* $g$*, then* $G \in \mathcal{G}(3, 6(1-g))$*.*

PROOF. Note that (ii) is a special case of (iii), and (iv) is a special case of (v). For (i), if $G_S = (S, E_S)$ is an induced subgraph of $G$, then $G_S$ cannot be $\Delta$-regular, and $|E_S| \le \frac{\Delta}{2}|S| - 1$. For (iii) and (v), the graph properties of having tree-width at most $t$, or genus at most $g$, are hereditary. Also, if $|V| = n$, a graph of tree-width $t$ has at most $tn - t(t+1)/2$ edges (see, e.g., [2], Theorem 1, Theorem 34), and a graph of genus $g$ at most $3n - 6(1-g)$ edges (see, e.g., [7], Theorem 7.5, Corollary 7.9). □

REMARK 8. In (i)–(iv) of Lemma 39, we have $b > 0$, but observe that in (v) we have $b > 0$ if $g = 0$ (planar), $b = 0$ if $g = 1$ (toroidal) and $b < 0$ if $g > 1$.

COROLLARY 40. *If a graph* $G = (V, E)$ *on* $n$ *vertices is:*

(i) *a nonregular connected graph with maximum degree* $\Delta$*, then* $\kappa(G) \le \frac{\Delta}{2} - \frac{1}{n}$*;*
(ii) *a forest, then* $\kappa(G) \le 1 - \frac{1}{n}$*;*
(iii) *a graph of tree-width* $t$*, then* $\kappa(G) \le t - \frac{t(t+1)}{2n}$*;*
(iv) *a planar graph, then* $\kappa(G) \le 3 - \frac{6}{n}$*;*
(v) *a graph of genus* $g > 0$*, let* $k_g = \frac{7}{2} + \sqrt{12g + \frac{1}{4}}$*, then*

$$\kappa(G) \le \kappa_g = \max\{(\lfloor k_g \rfloor - 1)/2, 3 + 6(g-1)/\lceil k_g \rceil\}.$$

PROOF. Follows directly from Lemmas 38 and 39. □

REMARK 9. Suppose that $g$ is chosen so that $k_g$ is an integer. The bound in Corollary 40(v) gives $\kappa_g = (k_g - 1)/2$ (because $k_g$ is the point at which the two arguments to the maximum are equal). The bound says that for every graph $G$ with genus $g$, $\kappa(G) \le \kappa_g$. This bound is tight because there is a graph $G$ with density $\kappa(G) = \kappa_g$ and genus $g$. In particular, the complete graph $K_{k_g}$ has density $\kappa_g$. If $k_g \ge 3$, it also has genus $g$. The smallest genus of a surface in which it can be embedded is $\gamma = \lceil (k_g - 3)(k_g - 4)/12 \rceil$ (see, e.g., [7], Theorem 7.10). This is at least[1] $g$ since

$$\gamma \ge \frac{k_g^2 - 7k_g + 12}{12} = g,$$

so the genus of $G$ is $g$ as required. The bound in Corollary 40(v) may not be tight for those $g$ for which $k_g$ is not integral. However, the bound is not greatly in error. Consider any $g > 0$. The graph $G = K_{\lfloor k_g \rfloor}$ can be embedded in a surface of genus $g$ so it has genus $g$. Also, as noted above, $\kappa(G) = \frac{1}{2}(\lfloor k_g \rfloor - 1)$. If the bound is not tight for this $g$ and $G$ then

$$\begin{aligned} \kappa(G) &\le \kappa_g = 3 + \frac{6(g-1)}{\lceil k_g \rceil} \le 3 + \frac{6(g-1)}{k_g} \\ &= \frac{k_g - 1}{2} = \left(\frac{\lfloor k_g \rfloor - 1}{2}\right) + \left(\frac{k_g - \lfloor k_g \rfloor}{2}\right) \le \kappa(G) + \frac{1}{2}, \end{aligned} \tag{7}$$

so $\kappa_g$ cannot be too much bigger than $\kappa(G)$. It is easy to see that $\kappa_g \sim \sqrt{3g}$ for large $g$. For small $g$, a plot of the upper bound $\kappa_g$ on maximum density is shown in Figure 2, together with the lower bound $\frac{1}{2}(\lfloor k_g \rfloor - 1)$.

[1]In fact, $\gamma = g$, though we do not use this fact here.

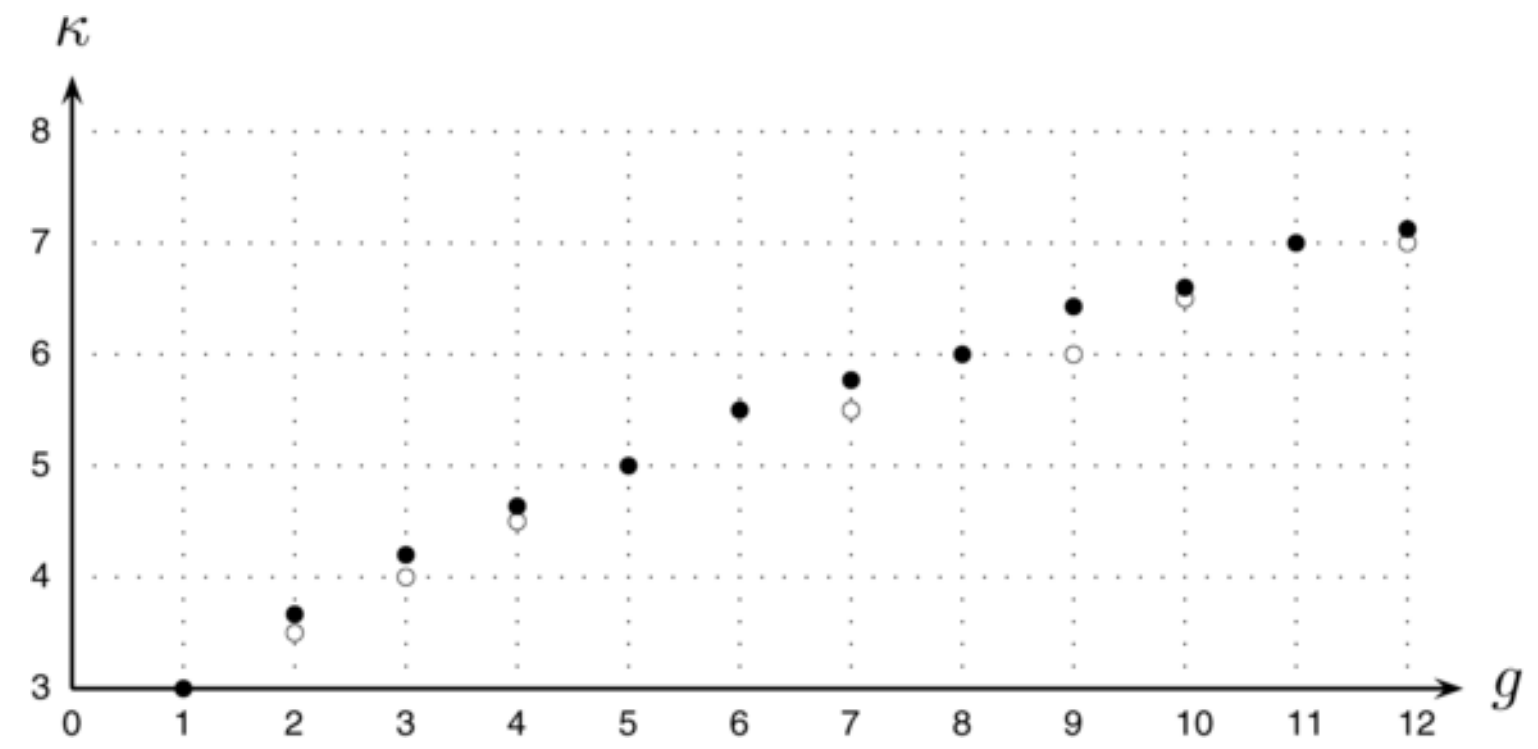

FIG. 2. *Upper and lower bounds on maximum density for small genus $g$.*

We now show that there exists a suitable $B$ for applying Corollary 37.

LEMMA 41. *Let $R \in \mathbb{M}_n^+$ be symmetric with maximum density $\kappa$ and let $\alpha = \|R\|_1$. Then there exists $B \in \mathbb{M}_n^+$ such that $R = B + B^{\mathrm{T}}$ and $\|B\|_1 = \kappa$, $\|B\|_\infty = \alpha - \kappa$.*

PROOF. It will be sufficient to show that $\|B\|_1 \le \kappa$, $\|B\|_\infty \le \alpha - \kappa$, since then we have

$$\begin{aligned} \alpha = \|R\|_1 &= \|B + B^{\mathrm{T}}\|_1 \le \|B\|_1 + \|B^{\mathrm{T}}\|_1 \\ &= \|B\|_1 + \|B\|_\infty \le \kappa + (\alpha - \kappa) = \alpha. \end{aligned} \tag{8}$$

First suppose $R$ is rational. Note that $\kappa$ is then also rational. Let $R' = R - D$, where $D = \operatorname{diag}(\varrho_{ii})$. Thus, for some large enough integer $N > 0$, $A(G) = NR'$ is the adjacency matrix of an undirected multigraph $G = (V, E)$ with $V = [n]$, $ND$ is a matrix of even integers, and $N\kappa$ is an integer. Thus, provided $B$ is eventually rescaled to $B/N$, we may assume these properties hold for $R'$, $D$ and $\kappa$. An *orientation* of $G$ is a directed multigraph $\vec{G} = (V, \vec{E})$ such that exactly one of $e^+ = (v, w)$, $e^- = (w, v)$ is in $\vec{E}$ for every $e = \{v, w\} \in E$. Clearly, $A(G) = A(\vec{G}) + A(\vec{G})^{\mathrm{T}}$, so we may take $B = A(\vec{G}) + \frac{1}{2}D$. Note that $\|B\|_1 = \max_{v\in V}(d_v^- + \frac{1}{2}\varrho_{vv})$ and $\|B\|_\infty = \max_{v\in V}(d_v^+ + \frac{1}{2}\varrho_{vv})$. We now apply the following (slightly restated) theorem of Frank and Gyárfás [20].

Theorem 42 (Frank and Gyárfás [20]). *Suppose $\ell_v \le u_v$ for all $v \in V$ in an undirected multigraph $G = (V, E)$. Then $G$ has an orientation $\vec{G}$ satisfying $\ell_v \le d_v^- \le u_v$ if and only if, for all $S \subseteq V$, we have $|E_S| \le \max\{\sum_{v\in S} u_v, \sum_{v\in S}(d_v - \ell_v)\}$.*

We will take $u_v = \kappa - \frac{1}{2}\varrho_{vv}$, $\ell_v = d_v + \kappa - \alpha + \frac{1}{2}\varrho_{vv}$. Then $\ell_v \le u_v$, since $d_v \le (\alpha - \varrho_{vv})$, and $(d_v - \ell_v) \ge u_v$, since $\alpha \ge 2\kappa$. The conditions of Theorem 42 are satisfied, since for all $S \subseteq V$,

$$\begin{aligned} |E_S| &= \tfrac{1}{2}\sum_{v\in S}\sum_{w\in S}\varrho_{vw} - \tfrac{1}{2}\sum_{v\in S}\varrho_{vv} \\ &\le \kappa|S| - \tfrac{1}{2}\sum_{v\in S}\varrho_{vv} = \sum_{v\in S} u_v \;\le\; \sum_{v\in S}(d_v - \ell_v). \end{aligned}$$

The result now follows for rational $R$, since we have

$$\begin{aligned} \|B\|_1 &= \max_{v\in V}(d_v^- + \tfrac{1}{2}\varrho_{vv}) \le \kappa, \\ \|B\|_\infty &= \max_{v\in V}(d_v^+ + \tfrac{1}{2}\varrho_{vv}) \le \max_{v\in V}(d_v - \ell_v + \tfrac{1}{2}\varrho_{vv}) = \alpha - \kappa. \end{aligned}$$

If $R$ is irrational, standard continuity arguments now give the conclusion. □

Remark 10. The use of Theorem 42 in the proof can be replaced by an application of the max-flow min-cut theorem, as in [21], but Theorem 42 seems more easily applicable here.

We can show that Lemma 41 is best possible, in the following sense.

Lemma 43. *Let $R \in \mathbb{M}_n^+$ be symmetric with maximum density $\kappa$ and let $\alpha = \|R\|_1$. If $R = B + B^{\mathrm{T}}$ for any $B \in \mathbb{M}_n$, then $\|B\|_1 \ge \kappa$ and $\|B\|_\infty \ge \alpha - \|B\|_1$.*

Proof. Let $I$ be any set achieving the maximum density of $R$. Then

$$\begin{aligned} 2|I|\kappa &= \sum_{i,j\in I}\varrho_{ij} \;\le\; \sum_{i,j\in I}(|B_{ij}| + |B_{ji}|) \\ &= 2\sum_{i,j\in I}|B_{ij}| \;\le\; 2\sum_{j\in I}\sum_{i\in[n]}|B_{ij}| \;\le\; 2|I|\|B\|_1, \end{aligned}$$

so $\|B\|_1 \geq \kappa$. The second assertion follows from (8). $\square$

THEOREM 44. *If $R \in \mathbb{M}_n^+$ is a symmetric matrix with maximum density $\kappa$ and $\alpha = \|R\|_1$, then $\lambda(R) \leq 2\sqrt{\kappa(\alpha - \kappa)}$.*

PROOF. Follows directly from Corollary 37 and Lemma 41. $\square$

REMARK 11. Since $\kappa(\alpha - \kappa)$ is increasing for $\kappa \leq \alpha/2$, an upper bound $\kappa'$ can be used, as long as we ensure that $\kappa' \leq \alpha/2$.

REMARK 12. We can adapt this for asymmetric $R$ by considering the "symmetrization" $\frac{1}{2}(R + R^{\mathrm{T}})$. Note that $\kappa(R) = \kappa(\frac{1}{2}(R + R^{\mathrm{T}}))$. Let $\tilde{\alpha}(R) = \|\frac{1}{2}(R + R^{\mathrm{T}})\|_1 \leq \frac{1}{2}(\|R\|_1 + \|R\|_\infty)$. We also have $\lambda(R) \leq \nu(R) = \nu(\frac{1}{2}(R + R^{\mathrm{T}})) = \lambda(\frac{1}{2}(R + R^{\mathrm{T}}))$. Then $\lambda(R) \leq 2\sqrt{\kappa(\tilde{\alpha} - \kappa)}$.

The following application, used together with Lemma 2, strengthens [14], Theorem 15.

THEOREM 45. *Suppose $R$ is a symmetric and irreducible dependency matrix with row sums at most 1, and suppose $0 < \gamma \leq \min_{i,j \in [n]}\{\varrho_{ij} : \varrho_{ij} > 0\}$. If there is any row with sum at most $1 - \gamma$, then $\lambda(R) \leq \sqrt{1 - \gamma^2/n^2} \leq 1 - \gamma^2/2n^2$.*

PROOF. Since $R$ is irreducible, for any $I \subset [n]$, $\sum_{i,j \in I} \varrho_{ij} \leq |I| - \gamma$. This also holds for $I = [n]$ by assumption. Thus, $\kappa \leq \frac{1}{2} - \frac{\gamma}{2n}$. Since $\|R\|_1 \leq 1$, we have $\lambda(R) \leq 2\sqrt{(\frac{1}{2} - \frac{\gamma}{2n})(\frac{1}{2} + \frac{\gamma}{2n})} = \sqrt{1 - \gamma^2/n^2}$. The final inequality is easily verified. $\square$

We can also apply Theorem 44 straightforwardly to (simple) graphs.

COROLLARY 46. *If $G$ has maximum density $\kappa$ and maximum degree $\Delta$, then $\lambda(G) \leq 2\sqrt{\kappa(\Delta - \kappa)}$.*

PROOF. In Theorem 44, we have $\alpha = \Delta$. $\square$

THEOREM 47. *If $G = (V, E) \in \mathcal{G}(a, b)$, with $b \geq 0$, $\Delta \geq 2a$ and $|V| = n$, then*

$$\lambda(G) \leq 2\sqrt{\left(a - \frac{b}{n}\right)\left(\Delta - a + \frac{b}{n}\right)}$$

$$\leq \begin{cases} \sqrt{a(\Delta - a)}\left(2 - \dfrac{b(\Delta - 2a)}{a(\Delta - a)n}\right), & \text{if } \Delta > 2a; \\ a\left(2 - \dfrac{b^2}{a^2 n^2}\right), & \text{if } \Delta = 2a. \end{cases}$$

PROOF. The first inequality follows directly from Lemma 38 and Corollary 46. Note that the condition $\Delta \geq 2a - 2b/n$ is required in view of Remark 11. For the second, squaring gives

$$4a(\Delta - a) - \frac{4b(\Delta - 2a)}{n} - \frac{4b^2}{n^2} \leq 4a(\Delta - a) - \frac{4b(\Delta - 2a)}{n} + \frac{b^2(\Delta - 2a)^2}{a(\Delta - a)n^2},$$

which holds for all $b$ and $\Delta \geq 2a$. When $\Delta = 2a$, using $\sqrt{1-x} \leq 1 - \frac{1}{2}x$,

$$\lambda(G) \leq 2\sqrt{a^2 - \frac{b^2}{n^2}} = 2a\sqrt{1 - \frac{b^2}{a^2n^2}}$$
$$\leq 2a\left(1 - \frac{b^2}{2a^2n^2}\right) = a\left(2 - \frac{b^2}{a^2n^2}\right). \qquad \square$$

THEOREM 48. *If $G = (V, E) \in \mathcal{G}(a, b)$, with $b \leq 0$ and $|V| = n$, let $k^* = a + \frac{1}{2} + \sqrt{(a + \frac{1}{2})^2 - 2b}$ and $\kappa^* = \max\{(\lfloor k^* \rfloor - 1)/2, a - b/\lceil k^* \rceil\}$. Then, if $\Delta \geq 2\kappa^*$, $\lambda(G) \leq \sqrt{\kappa^*(\Delta - \kappa^*)}$.*

PROOF. This follows immediately from Lemma 38, Theorem 44 and Remark 11. $\square$

We can apply this to the examples from Lemma 39.

COROLLARY 49. *If $G = (V, E)$, with maximum degree $\Delta$ and $|V| = n$, is:*

(i) *A connected nonregular graph, then*
$$\lambda(G) \leq \sqrt{\Delta^2 - \frac{4}{n^2}} < \Delta - \frac{2}{\Delta n^2}.$$

(ii) *A tree with $\Delta \geq 2$, then*
$$\lambda(G) \leq 2\sqrt{\left(1 - \frac{1}{n}\right)\left(\Delta - 1 + \frac{1}{n}\right)} < \sqrt{\Delta - 1}\left(2 - \frac{\Delta - 2}{(\Delta - 1)n}\right).$$

*If $\Delta = 2$, then $\lambda(G) < 2 - 1/n^2$, and if $\Delta < 2$, then $\lambda(G) = \Delta$.*

(iii) *A graph with tree-width at most $t$ and $\Delta \geq 2t$, then*
$$\lambda(G) \leq 2\sqrt{\left(t - \frac{t(t+1)}{2n}\right)\left(\Delta - t + \frac{t(t+1)}{2n}\right)}$$
$$< \sqrt{t(\Delta - t)}\left(2 - \frac{(t+1)(\Delta - 2t)}{2(\Delta - t)n}\right).$$

*If $\Delta = 2t$, then $\lambda(G) < 2t - t(t+1)^2/4n^2$.*

(iv) *A planar graph with* $\Delta \geq 6$*, then*

$$\lambda(G) \leq 2\sqrt{(3-6/n)(\Delta-3+6/n)} < 2\sqrt{3(\Delta-3)}\left(1-\frac{\Delta-6}{(\Delta-3)n}\right).$$

*If* $\Delta = 6$*,* $\lambda(G) \leq 6 - 12/n^2$*. If* $\Delta \leq 5$*,* $\lambda(G) \leq \Delta$ *is best possible.*

(v) *A graph of genus* $g > 0$*, let* $k_g = \frac{7}{2} + \sqrt{12g + \frac{1}{4}}$ *and* $\kappa_g = \max\{(\lfloor k_g \rfloor - 1)/2, 3 + 6(g-1)/\lceil k_g \rceil\}$*. If* $\Delta \geq 2\kappa_g$*, then*

$$\lambda(G) \leq \sqrt{\kappa_g(\Delta - \kappa_g)}.$$

Proof. Using Lemma 39, these follow using Theorem 47 and Theorem 48 with:

(i) $a = \Delta/2,\ b = 1$ and $\Delta = 2a$;
(ii) $a = 1,\ b = 1$, if $\Delta > 2$. If $\Delta = 2$, the result follows from the $\Delta = 2a$ case. $\Delta = 1$, $G$ is a single edge and, if $\Delta = 0$, an isolated vertex;
(iii) $a = t,\ b = t(t+1)/2$;
(iv) $a = 3,\ b = 6$. If $\Delta \leq 5$, regular planar graphs with degree $\Delta$ exist, and we use Lemma 35;
(v) $a = 3,\ b = -6(g-1)$. □

Remark 13. If $G$ is a disconnected graph, the component having the largest eigenvalue determines $\lambda(G)$, using Lemma 8. This can be applied to a forest.

Remark 14. Corollary 49(i) improves on a result of Stevanović [37], who showed that

$$\lambda(G) \ < \ \Delta - \frac{1}{2n(n\Delta-1)\Delta^2}.$$

This was improved by Zhang [40] to (approximately) $\Delta - \frac{1}{2}(\Delta n)^{-2}$, which is still inferior to (i). But recently the bound has been improved further by Cioabă, Gregory and Nikiforov [6], who showed

$$\lambda(G) \ < \ \Delta - \frac{1}{n(\mathcal{D}+1)},$$

where $\mathcal{D}$ is the diameter of $G$. This gives $\lambda(G) \leq \Delta - 1/n^2$ even in the worst case, which significantly improves on (i). However, Corollary 49 is an easy consequence of the general Corollary 46, whereas [6] uses a calculation carefully tailored for this application.

Remark 15. When $G$ is a degree-bounded forest, Corollary 49(ii) strengthens another result of Stevanović [36], who showed $\lambda(G) < 2\sqrt{\Delta - 1}$.

Remark 16. When $G$ is a planar graph, Theorem 47(iv) improves a result of Hayes [22].

We can now apply these results to the mixing of Glauber dynamics for proper colorings in the classes of sparse graphs $\mathcal{G}(a,b)$.

THEOREM 50. *Let $G=(V,E)\in\mathcal{G}(a,b)$, with $b>0$, have maximum degree $\Delta\geq 2a$, where $|V|=n$. Let $\psi=2\sqrt{a(\Delta-a)}$, $\phi=\Delta-2a$ and $\mu=\psi/(q-\Delta)$. Then, if:*

(i) *$q>\Delta+\psi$, the random update and systematic scan Glauber dynamics mix in time*

$$\tau_{\mathrm{r}}(\varepsilon)\ \leq\ (1-\mu)^{-1}\,n\ln(n/\varepsilon),\qquad \hat{\tau}_{\mathrm{s}}(\varepsilon)\sim(1-\tfrac{1}{2}\mu)(1-\mu)^{-1}\ln(n/\varepsilon).$$

(ii) *$q=\Delta+\psi$ and $\phi>0$, the random update and systematic scan Glauber dynamics mix in time*

$$\tau_{\mathrm{r}}(\varepsilon)\ \leq\ (\psi^2/2b\phi)n^2\ln(n/\varepsilon),\qquad \hat{\tau}_{\mathrm{s}}(\varepsilon)\sim(\psi^2/2b\phi)n\ln(n/\varepsilon).$$

(iii) *$q=\Delta+\psi$ and $\phi=0$, the random update and systematic scan Glauber dynamics mix in time*

$$\tau_{\mathrm{r}}(\varepsilon)\leq 2(a/b)^2n^3\ln(n/\varepsilon),\qquad \hat{\tau}_{\mathrm{s}}(\varepsilon)\sim 3(a/b)^2n^2\ln(n/\varepsilon).$$

PROOF. Recall from the beginning of Section 4 that $\lambda(R)\leq\lambda(G)/(q-\Delta)$ where $\lambda(G)$ denotes $\lambda(A(G))$. Note also that, if $\psi$ is not an integer, then $q-\Delta-\psi=\Omega(1)$. By Theorem 47, for (i) we have $\|R\|_2=\lambda(R)\leq\lambda(G)/(q-\Delta)\leq\psi/(q-\Delta)=\mu<1$. For (ii), we have $\lambda(R)\leq 1-(2b\phi/\psi^2n)$, and for (iii), $\lambda(R)\leq 1-(b^2/2a^2n^2)$. The conclusions for $\tau_{\mathrm{r}}(\varepsilon)$ follow from Lemma 17, and those for $\tau_{\mathrm{s}}(\varepsilon)$ from Lemma 3. For (ii) and (iii), factors of $\frac{1}{2}$ arise in Lemma 3 since $\lambda\sim 1$, but additional factors (2 and 3, resp.) come from the log term. □

THEOREM 51. *If $G=(V,E)\in\mathcal{G}(a,b)$ with $b\leq 0$, let $k^*=a+\frac{1}{2}+\sqrt{(a+\frac{1}{2})^2-2b}$ and $\kappa^*=\max\{(\lfloor k^*\rfloor-1)/2,a-b/\lceil k^*\rceil\}$. If $\Delta>2\kappa^*$, let $\psi=\sqrt{\kappa^*(q-\kappa^*)}$ and $\mu=\psi/(q-\Delta)$. Then, if $q>\Delta+\psi$,*

$$\tau_r(\varepsilon)\ \leq\ (1-\mu)^{-1}\,n\ln(n/\varepsilon),\qquad \hat{\tau}_{\mathrm{s}}(\varepsilon)\sim(1-\tfrac{1}{2}\mu)(1-\mu)^{-1}\ln(n/\varepsilon).$$

PROOF. From Theorem 48, we have

$$\|R\|_2=\lambda(R)\leq\frac{\lambda(G)}{q-\Delta}\leq\frac{\psi}{q-\Delta}=\mu<1.$$

The conclusions for $\tau_{\mathrm{r}}(\varepsilon)$ now follow from Lemmas 14 and 17, and those for $\hat{\tau}_{\mathrm{s}}(\varepsilon)$ from Lemma 3. □

COROLLARY 52. *If $G=(V,E)$, with $|V|=n$ and maximum degree $\Delta$, is:*

(i) *a nonregular connected graph and $q=2\Delta$, then*

$$\tau_{\mathrm{r}}(\varepsilon)\ \leq\ \tfrac{1}{2}\Delta^2n^3\ln(n/\varepsilon),\qquad \hat{\tau}_{\mathrm{s}}(\varepsilon)\ \sim\ \tfrac{3}{4}\Delta^2n^2\ln(n/\varepsilon).$$

(ii) *A graph with tree-width $t$ and $\Delta \geq 2t$, let $\psi = 2\sqrt{t(\Delta - t)}$. Then*

$$\tau_{\mathrm{r}}(\varepsilon) \leq \begin{cases} (q-\Delta)(q-\Delta-\psi)^{-1} n \ln(n/\varepsilon), & \text{if } q > \Delta + \psi; \\ \psi^2(t(t+1)(\Delta - 2t))^{-1} n^2 \ln(n/\varepsilon), & \text{if } q = \Delta + \psi \text{ and } \Delta > 2t; \\ 8(t+1)^{-2} n^3 \ln(n/\varepsilon), & \text{if } q = \Delta + \psi \text{ and } \Delta = 2t. \end{cases}$$

$$\hat{\tau}_{\mathrm{s}}(\varepsilon) \sim \begin{cases} (q-\Delta-\frac{1}{2}\psi)(q-\Delta-\psi)^{-1} \ln(n/\varepsilon), & \text{if } q > \Delta + \psi; \\ \psi^2(t(t+1)(\Delta - 2t))^{-1} n \ln(n/\varepsilon), & \text{if } q = \Delta + \psi \text{ and } \Delta > 2t; \\ 12(t+1)^{-2} n^2 \ln(n/\varepsilon), & \text{if } q = \Delta + \psi \text{ and } \Delta = 2t. \end{cases}$$

(iii) *A planar graph and $\Delta \geq 6$, let $\psi = 2\sqrt{3(\Delta - 3)}$. Then*

$$\tau_{\mathrm{r}}(\varepsilon) \leq \begin{cases} (q-\Delta)(q-\Delta-\psi)^{-1} n \ln(n/\varepsilon), & \text{if } q > \Delta + \psi; \\ \psi^2(12(\Delta - 6))^{-1} n^2 \ln(n/\varepsilon), & \text{if } q = \Delta + \psi \text{ and } \Delta > 6; \\ \frac{1}{2} n^3 \ln(n/\varepsilon), & \text{if } q = \Delta + \psi \text{ and } \Delta = 6. \end{cases}$$

$$\hat{\tau}_{\mathrm{s}}(\varepsilon) \sim \begin{cases} (q-\Delta-\frac{1}{2}\psi)(q-\Delta-\psi)^{-1} \ln(n/\varepsilon), & \text{if } q > \Delta + \psi; \\ \psi^2(12(\Delta - 6))^{-1} n \ln(n/\varepsilon), & \text{if } q = \Delta + \psi \text{ and } \Delta > 6; \\ \frac{3}{4} n^2 \ln(n/\varepsilon), & \text{if } q = \Delta + \psi \text{ and } \Delta = 6. \end{cases}$$

(iv) *A graph of genus $g > 0$, let $k_g = \frac{7}{2} + \sqrt{12g + \frac{1}{4}}$, $\kappa_g = \max\{(\lfloor k_g \rfloor - 1)/2, 3 + 6(g-1)/\lceil k_g \rceil\}$ and $\psi = \sqrt{\kappa_g(\Delta - \kappa_g)}$. If $\Delta > 2\kappa_g$ and $q > \Delta + \psi$, then*

$$\tau_{\mathrm{r}}(\varepsilon) \leq (q-\Delta)(q-\Delta-\psi)^{-1} n \ln(n/\varepsilon),$$
$$\hat{\tau}_{\mathrm{s}}(\varepsilon) \sim (q-\Delta-\tfrac{1}{2}\psi)(q-\Delta-\psi)^{-1} \ln(n/\varepsilon).$$

PROOF. This follows directly from Lemma 39 and Theorems 50 and 51. □

REMARK 17. Corollary 52(i) bounds the mixing time of heat bath Glauber dynamics for sampling proper $q$-colorings of a *nonregular* graph $G$ with maximum degree $\Delta$ when $q = 2\Delta$. (We can bound the mixing time for a disconnected graph $G$ by considering the components.) It is also possible to extend the mixing time result for nonregular graphs to regular graphs using the decomposition method of Martin and Randall [26]. See [14], Section 5, for details about how to do this. The use of our Corollary 52(i) improves Theorem 5 of [14] by a factor of $n$.

M. Dyer
School of Computing
University of Leeds
Leeds LS2 9JT
United Kingdom

L. A. Goldberg
Department of Computer Science
University of Liverpool
Liverpool L69 3BX
United Kingdom
E-mail: L.A.Goldberg@liverpool.ac.uk

M. Jerrum
School of Mathematical Sciences
Queen Mary, University of London
Mile End Road, London E1 4NS
United Kingdom